\documentclass{fic-l}
\usepackage{amsmath,amssymb,amsfonts,multirow,epsfig,xspace,verbatim}
\usepackage[round]{natbib}

\newcommand{\sref}[1]{\S~\ref{sec:#1}}

\newcommand{\fref}[1]{Figure~\ref{fig:#1}}

\newcommand{\R}{{\mathbb R}}
\newcommand{\Tmix}{T_{\operatorname{mix}}}
\newcommand{\Tsep}{T_{\operatorname{sep}}}
\newenvironment{coupling}{\begin{list}{\hspace*{0em}}{\renewcommand{\parsep}{-2pt}}}{\end{list}}
\newenvironment{code}{\begin{list}{\hspace*{0em}}{\renewcommand{\parsep}{-2pt}}}{\end{list}}

\newcommand{\Repeat}{{\rm repeat\ }}
\newcommand{\Until}{{\rm until\ }}
\newcommand{\For}{{\rm for\ }}
\newcommand{\To}{{\rm to\:\,}}
\newcommand{\DownTo}{{\rm downto\:\,}}
\newcommand{\If}{{\rm if\ }}

\newcommand{\Output}{{\rm output\ }}

\newcommand{\s}{\hspace*{3ex}}
\newcommand{\sspace}{$\langle\text{state space}\rangle$}

\newcommand{\epf}{\textit{ex post facto}\xspace}
\newcommand{\Epf}{\textit{Ex post facto}\xspace}
\newcommand{\PDF}{\operatorname{PDF}}
\newcommand{\LIPDF}{\operatorname{LeftInversePDF}}
\newcommand{\RIPDF}{\operatorname{RightInversePDF}}
\newcommand{\Mode}{\operatorname{Mode}}
\newcommand{\Uniform}{\operatorname{Uniform}}
\newcommand{\Exponential}{\operatorname{Exponential}}
\newcommand{\Normal}{\operatorname{Normal}}
\newcommand{\GammaD}{\operatorname{Gamma}}

\newcommand{\fslrx}   {$\displaystyle f_{L,R,X}(s) = \left\lfloor \frac{s+R-X}{R-L}\right\rfloor (R-L) +X$}
\newcommand{\fslrxepf}{$\displaystyle f_{L,R,X_0}(s) = \left\lfloor \frac{s+R-X_0}{R-L}\right\rfloor (R-L) +X_0$}

\newtheorem{theorem}{Theorem}[section]

\theoremstyle{definition}

\theoremstyle{remark}

\numberwithin{equation}{section}

\newtheorem{claim}[theorem]{Claim}

\makeatletter
\newcommand\@dotsep{4.5}
\renewcommand*\l@section[2]{%
  \ifnum \c@tocdepth >\z@
    \addpenalty\@secpenalty
    \addvspace{1.0em \@plus\p@}%
    \setlength\@tempdima{1.5em}%
    \begingroup
      \parindent \z@ \rightskip \@pnumwidth
      \parfillskip -\@pnumwidth
      \leavevmode \bfseries
      \advance\leftskip\@tempdima
      \hskip -\leftskip
      #1\nobreak\hfil \nobreak\hb@xt@\@pnumwidth{\hss #2}\par
    \endgroup
  \fi}
\renewcommand*\l@subsection{\@dottedtocline{2}{1.5em}{2.3em}}
\renewcommand*\l@subsubsection{\@dottedtocline{3}{3.8em}{3.2em}}
\def\@serieslogo{}
\makeatother
\begin{document}

\title{Layered Multishift Coupling\\
 for use in Perfect Sampling Algorithms\\
 (with a primer on CFTP)}

\author{David Bruce Wilson}
\address{Microsoft Research\\One Microsoft Way\\Redmond, WA 98052}
\email{dbwilson@alum.mit.edu}
\thanks{The research that led to this article was done in part while at the Institute for Advanced Study (supported by the NSF), and in part while at Microsoft.}

\renewcommand{\subjclassname}{\textup{2000} Mathematics Subject Classification}
\subjclass{65C40}
\date{}

\begin{abstract}
In this article we describe a new coupling technique which is useful
in a variety of perfect sampling algorithms.  A multishift coupler
generates a random function $f()$ so that for each $x\in{\mathbb{R}}$,
$f(x)-x$ is governed by the same fixed probability distribution, such
as a normal distribution.  We develop the class of layered multishift
couplers, which are simple and have several useful properties.  For
the standard normal distribution, for instance, the layered multishift
coupler generates an $f()$ which (surprisingly) maps an interval of
length $\ell$ to fewer than $2+\ell/2.35$ points --- useful in
applications which perform computations on each such image point.  The
layered multishift coupler improves and simplifies algorithms for
generating perfectly random samples from several distributions,
including the autogamma distribution, posterior
distributions for Bayesian inference, and the steady state
distribution for certain storage systems.  We also use the layered
multishift coupler to develop a Markov-chain based perfect sampling
algorithm for the autonormal distribution.

At the request of the organizers, we begin by giving a primer on CFTP
(coupling from the past); CFTP and Fill's algorithm are the two
predominant techniques for generating perfectly random samples using
coupled Markov chains.
\end{abstract}

\maketitle
\markright{Layered Multishift Coupling for use in Perfect Sampling Algorithms}

\addtolength{\textheight}{3\baselineskip}
\newpage

\tableofcontents
\def\contentsname{Contents}
\addtolength{\textheight}{-3\baselineskip}
\newpage
\makeatletter
\renewcommand*\l@section{\@tocline{1}{0pt}{1pc}{}{}}
\renewcommand*\l@subsection{\@tocline{2}{0pt}{1pc}{5pc}{}}
\renewcommand*\l@subsubsection{\@tocline{3}{0pt}{1pc}{7pc}{}}
\makeatother

\section{Primer on Coupling from the Past}

\subsection{Markov chain Monte Carlo}

For many applications in statistical physics, computer science, and
Bayesian inference, it is very useful to generate random structures
according to some pre-specified distribution.  Sometimes there is a
direct random generation method, such as with percolation, random
permutations, or Gaussian random variables.  But often the state
spaces are more complicated and there is no known direct sampling
method, as is the case for random independent sets of a graph, random
linear extensions of a partially ordered set, or random contingency
tables.  To sample from state spaces such as these, people typically
rely upon Markov chains.  There is some natural randomizing operation,
which given an input state, produces a randomly modified output state.
If the input state is already distributed according to the desired
distribution, then so is the output state.  Under mild conditions, if
sufficiently many randomizing operations are performed, then the final
state will be distributed in approximately the desired distribution.

The computer which simulates the Markov chain doesn't have any
idea what ``sufficiently many'' means.  This may mean
one of the following.
\begin{itemize}
\item
 The computer keeps simulating the Markov chain forever.  This may be
 OK when doing mathematics, but it is not a practical approach to MCMC.
\item
 The human guesses how many steps must be enough.  The guess could be bad.
\item
 The human applies spectral analysis or other mathematical techniques
 to rigorously determine how many steps are enough.  Obtaining
 rigorous bounds has been an active area of research in the past
 decade, and there have been some notable successes.  More than one
 person has been elected to the National Academy of Sciences for work
 in this area.  But by and large this remains a hard problem, and many
 (if not most) Markov chains of practical interest have so far failed
 to succumb to rigorous analysis.
\item
 The human writes code to measure various autocorrelation functions,
 thereby allowing the computer to heuristically guess how many steps
 are enough.  This method is the workhorse for MCMC in physics and
 statistics.  In the absence of a better option, it gets the job done.
 But no matter how good the heuristics are, one can never be
 completely sure that they did the job correctly, and the allegedly
 random samples produced could be quite biased.
\item
 The computer (rigorously) figures out on its own how long to run the
 Markov chain.  When it is possible to do this, life is simplified for
 the experimenter. The two predominant methods for doing this are
 known as ``coupling from the past'' (CFTP)
 \citep{propp-wilson:exact-sampling} and Fill's algorithm
 \citep{fill:interruptible}.  The first of these methods is the topic
 of this primer; additional information on the second is given by
 \cite*{fill-machida-murdoch-rosenthal:fill}.
\end{itemize}

\subsection{Randomizing operations, pairwise couplings, and Markov chains}
\label{sec:randomizing}

Typically a randomizing operation is expressed as a deterministic
function $\phi$ that takes two inputs, the input state $X_t$ at time
$t$ and some intrinsic randomness $U_t$, and returns the modified
output state $X_{t+1}=\phi(X_t,U_t)$.  One can think of $\phi$ as
representing a piece of C code or Lisp code, and $U_t$ as representing
the output of the pseudorandom number generator.  It is assumed that
the $U_t$'s are mutually independent.  Conceptually it is convenient
to combine $\phi$ and $U_t$ into a single random function $f_t$
defined by $f_t(x) = \phi(x,U_t)$.  The randomizing operation is
assumed to preserve the desired distribution $\pi$ from which we wish
to sample: if $X_t$ is distributed according to $\pi$ and $U_t$ is
random, then $X_{t+1}$ is also distributed according to $\pi$.

Applying the randomizing operation to a given
state is equivalent to running the Markov chain one step from that
state.  There can be many different randomizing operations that are
consistent with a given Markov chain.

Just as a toy example, suppose that the state space consists of the
integers from $0$ to $n$, that $U_t$ is $+1$ or $-1$ with probability
$1/2$ each, and $\phi$ is defined by
$$\phi(x,u) = \begin{cases} x+u & 0\leq x+u \leq n\\ 0 & x+u<0\\ n & x+u>n\end{cases}.$$
Then it is easy to check that the only distribution $\pi$ preserved by
these randomizing operations is the uniform distribution on
$0,1,\ldots,n$.

A different randomizing operation $\psi$ might flip $n+1$ different
coins, and use the $x$th coin flip when computing
$\psi(x,u)=\phi(x,u_x)$.  While these two randomizing operations are
different, they give rise to the same Markov chain.  For CFTP, the
choice of randomizing operation is as important as the Markov chain
itself.

A pairwise coupling is a method for updating a pair of states so that
the evolution of either state by itself is described by the Markov
chain.  Randomizing operations are also sometimes called
``simultaneous couplings'' or ``grand couplings'', since the specified
how any group of states will evolve.  There are exceptions (see e.g.\
\sref{monotone}), but a large fraction of the pairwise couplings
encountered in practice extend naturally to grand couplings /
randomizing operations.  For the purposes of CFTP, we will principally
be interested in randomizing operations.

\subsection{CFTP: Sampling with and then without an oracle}
\label{sec:cftp}

Suppose that we have an oracle which returns perfectly random samples
distributed according to $\pi$.  Such oracle would make life easy for
someone doing Monte Carlo experiments, if not for the fact that it
charges \$20 for each sample requested of it.  Since we do not have an
unlimited budget, we would like to make use of our randomizing
operation, which is essentially free, and thereby reduce our
dependence upon the oracle.

To this end, consider experiment $A_T$ given below:
\begin{code}
\item   Pay \$20 to the oracle to draw $X_{-T}$ from $\pi$
\item   For $t=-T$ upto $-1$
\item\s   Compute $X_{t+1} := \phi(X_t,U_t)$
\item   Output $X_0$
\end{code}
Because the distribution $\pi$ is preserved by the randomizing
operation, by induction follows that the output state $X_0$ is
distributed exactly according to $\pi$.  Next consider experiment
$B_T$ given below:
\begin{code}
\item   Look at $U_{-T},U_{-T+1},\ldots,U_{-1}$
\item   If there is only one possible value for $X_0$
\item   Then
\item\s   Output $X_0$
\item   Else 
\item\s   Pay \$20 to the oracle to draw $X_{-T}$ from $\pi$
\item\s   For $t=-T$ upto $-1$
\item\s\s   Compute $X_{t+1} := \phi(X_t,U_t)$
\item\s   Output $X_0$
\end{code}
If we're lucky, experiment $B_T$ does not need to pay the oracle \$20.
Note that experiment $B_T$ always returns precisely the same answer
that experiment $A_T$ does, provided that the same random values $U_t$
are used.  From this we see that the output of experiment $B_T$ is
distributed exactly according to $\pi$, provided that the random
values $U_t$ used in the second part are the same values used in the
first part of the procedure.  Using fresh random values in the second
part of the procedure is a bad idea that would cause the output state
to be biased.  Also note that the running time of experiment $B_T$ is
a random variable which may be correlated with the output state.  Thus
repeatedly running, interrupting, and restarting the procedure is also
a bad idea that would introduce bias --- much better to let the
procedure finish running and return its answer.

We return to our toy example to see how one might conduct experiment
$B_T$ in practice.  \fref{oracle} shows two possible outcomes of
experiment $B_{10}$.

\begin{figure}[phtb]
\centerline{\epsfig{figure=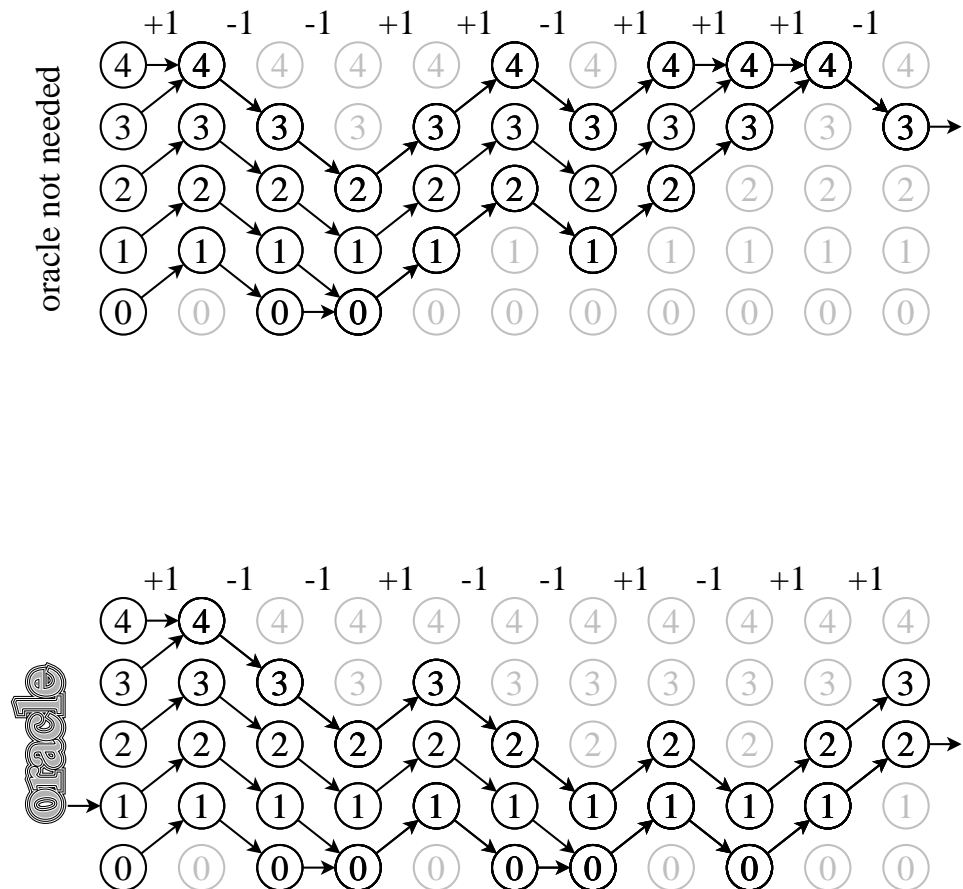,width=\textwidth}}
\caption{
The outcome of experiment $B_{10}$ using two different sequences of
the coins $U_{-10},\ldots,U_{-1}$.  The horizontal axis represents
time.  In the first example, the coins determine the final state
$X_0$, so there is no need to consult the oracle, and state $3$ is
output.  In the second example, the coins restrict the possible values
of $X_0$ but do not determine it.  The oracle is then consulted (for a
fee of \$20), and returns state $1$, which gets sent to state $2$ (the
output) by the coins.  In both cases we only had to track the top-most
and bottom-most trajectories (rather than all of them) to see if the
coins determine state $X_0$.  }
\label{fig:oracle}
\end{figure}

In the first example the random choices of $U_{-10},\ldots,U_{-1}$
were $+1$, $-1$, $-1$, $+1$, $+1$, $-1$, $+1$, $+1$, $+1$, $-1$.
Every possible starting value for state $X_{-10}$ was tried, and for
each of these starting values the final state $X_0$ was $3$.  Since
the given choices of $U_{-10},\ldots,U_{-1}$ determine $X_0$, the
oracle was not consulted, and the final output state was $3$.

This example also illustrates a convenient property of what are known
as ``monotone Markov chains''.  The state space comes equipped with a
partial order $\preceq$ with a biggest state $\hat{1}$ and a smallest
state $\hat{0}$ such that $\hat0\preceq x\preceq \hat1$ for all states
$x$.  In this toy example the partial order $\preceq$ is the usual
order $\leq$ on integers, $\hat1=4$, and $\hat0=0$.  A randomizing
operation on a partially ordered set is monotone if $x\preceq y$
implies $\phi(x,u) \preceq \phi(y,u)$.  For monotone Markov chains it
is particularly easy to test if the $U_t$'s determine $X_0$: apply the
randomizing operations starting from $X_{-T}=\hat0$ and then from
$X_{-T}=\hat1$.  If $f_{-1}(\cdots f_{-T}(\hat0)\cdots) =
f_{-1}(\cdots f_{-T}(\hat1)\cdots)$, then no matter what starting
value for $X_{-T}$ that the oracle would have selected, the final
value for $X_0$ is just $f_{-1}(\cdots f_{-T}(\hat1)\cdots)$.  So it
is only necessary to test two possible starting values rather than all
of them.

In the second example a different set of random choices of
$U_{-10},\ldots,U_{-1}$ is used.  In this case $U_{-10},\ldots,U_{-1}$
did not determine $X_0$, so \$20 was paid to the oracle, which then
assigned $X_{-10}=1$.  Applying the randomizing operations specified
by $U_{-10},\ldots,U_{-1}$ resulted in the final state $X_0=2$, which
was the output.

In order to reduce the chance that we have to resort to paying \$20 to
the oracle, we should pick a large value of $T$ when doing experiment
$B_T$, since that would increase the probability that
$U_{-T},\ldots,U_{-1}$ determine $X_0$.  But if we pick an excessively
large value of $T$, we'd rather not spend time looking at all the
$U_t$'s if the last several of them by themselves determine $X_0$.  We
could look at the last several $U_t$'s and see if they determine
$X_0$.  If not, we can continue to look at progressively more of the
$U_t$'s to see if they determine $X_0$.  If we are unlucky and
$U_{-T},\ldots,U_{-1}$ fail to determine $X_0$, then we resort to
paying \$20 to the oracle.  This strategy is expressed more formally
as experiment $C_T$ given below.  It is evident that experiment $C_T$
and experiment $A_T$ will return the same answer provided that they
use the same values of $U_{-T},\ldots,U_{-1}$ and the oracle return
the same sample (if asked to do so).  Thus experiment $C_T$ returns a
random sample distributed exactly according to $\pi$, assuming of
course that each time it looks at a given random variable $U_t$, it
sees the same value.  For this reason people often stress the
importance of ``re-using the same random coins''.
\begin{code}
\item   If $U_{-1}$ determines $X_0$
\item   Then Output $X_0$
\item   Else If $U_{-2},U_{-1}$ determine $X_0$
\item   Then Output $X_0$
\item   Else If $U_{-4},U{-3},U_{-2},U_{-1}$ determine $X_0$
\item   Then Output $X_0$
\item   \hspace{3.5em} $\vdots$
\item   Else If $U_{-2^{\lceil\log_2 T\rceil-1}},\ldots,U_{-1}$ determine $X_0$
\item   Then Output $X_0$
\item   Else If $U_{-T},\ldots,U_{-1}$ determine $X_0$
\item   Then Output $X_0$
\item   Else
\item\s   Pay \$20 to the oracle to draw $X_{-T}$ from $\pi$
\item\s   For $t=-T$ upto $-1$
\item\s\s   Compute $X_{t+1} := \phi(X_t,U_t)$
\item\s   Output $X_0$
\end{code}

Coupling from the past is experiment $C_\infty$, which is defined by
$$\text{CFTP} = \text{Experiment $C_\infty$}
              = \lim_{T\rightarrow\infty}\text{Experiment $C_T$}.$$
Since all experiments $C_T$ (for large enough $T$) start out by doing
the same thing, taking this sort of limit makes sense.  Experiment
$C_\infty$, which is re-expressed below, has the convenient property
that it never consults the oracle.  CFTP is also illustrated in \fref{mcftp}.
\begin{code}
\item   $T := 1$
\item   While $U_{-T},\ldots,U_{-1}$ do not determine $X_0$
\item\s   $T := 2*T$
\item   Output $X_0$
\end{code}

\begin{figure}[phtb]
\centerline{\epsfig{figure=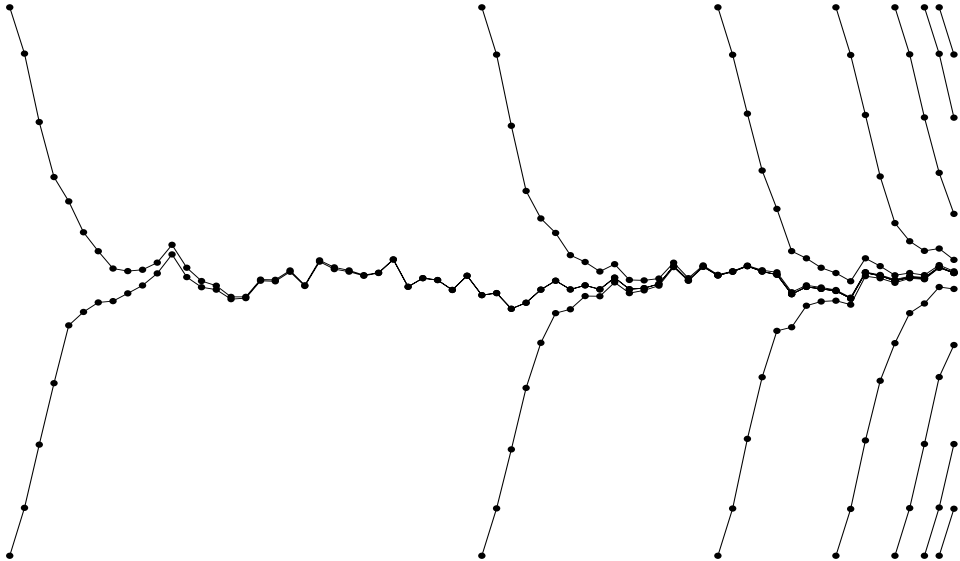,width=\textwidth,height=3in}}
\caption{
Illustration of CFTP in the monotone setting.  Shown are the heights
of the upper and lower trajectories started at various starting times
in the past.  When a given epoch is revisited later by the algorithm,
it uses the same randomizing operation.
}
\label{fig:mcftp}
\end{figure}

\subsection{Questions and answers}

When explaining CFTP to an audience, there are invariably many questions.  Included below are some of these questions together with their answers.

\newcommand{\q}{\smallskip\noindent\textbf{Q}: }
\newcommand{\qa}[1]{\smallskip\noindent\textbf{Q}#1: }
\newcommand{\ans}{

\noindent\textbf{A}: }
\newcommand{\ansa}[1]{

\noindent\textbf{A#1}: }

\q If we always end up in state $3$ no matter where we start, then how is
that a random sample?
\ans 
 For a given particular sequence of coin flips, every possible
 starting state ends up in state $3$.  But for a different (random)
 sequence of coin flips, every possible starting state may end up in a
 different final state at time $0$.

\q What if we just run the Markov chain forward until coalescence?
\ans Consider the toy example of the previous section.  Coalescence only occurs at the states $0$ and $n$, so the result would be very far from being distributed according to $\pi$.

\q What happens if we use fresh coins rather than re-using the same coins?
\ans
 Consider the toy example of the previous section, and set $n=2$ (so
 the states are $0,1,2$).  Then one can check that the probability of
 outputting state $1$ has a binary expansion of
 $0.0010010101001010101010101001\ldots$, which is neither rational nor
 very close to $1/3$.

\q Rather than doubling back in time, what if we double forward in time, and stop the Markov chain at the first power of $2$ greater than or equal to the coalescence time?
\ans
 As before, set $n=2$ in our toy example.  One can check that the
 probability of outputting state $1$ is $1/6$ rather than $1/3$.

\q What if we instead ...
\ans
 Enough already!  There do exist other ways to generate perfectly
 random samples, but the only obvious change that can be made to the
 CFTP algorithm without breaking it is changing the sequence of
 starting times in the past from powers of $2$ to some other sequences
 integers.

\q If coupling forwards in time doesn't work, then why does going backwards in time work?
\ans
 If you did not like the first explanation, then another way to look
 at it is that a virtual Markov chain has been running for all time,
 and so today the state is random.  If we can (with probability 1)
 figure out today's state by looking at some of the recent randomizing
 operations, then we have a random state.

\q If going backwards in time works, then why doesn't going forwards in time work?
\ansa{1}
 There is a way to obtain perfectly random samples by running forwards
 in time \citep{wilson:rocftp}, but none of the obvious variations
 described above work.
\ansa{2}
 CFTP determines the state at a deterministic time, any one of which
 is distributed according to $\pi$.  The variations suggested above
 determine the state at a random time, where that time is correlated
 with the moves of the Markov chain in a complicated way, and these
 correlations mess things up.

\q Why did you step back by powers of 2, when any other sequence would have worked?
\ans
For efficiency reasons.  Let $T^*$ be the best time in the past at
which to start, i.e.\ the smallest integer for which starting at time
$-T^*$ leads to coalescence.  By stepping back by powers of $2$, the
total number of Markov chain steps simulated is never larger than $4
T^*$ (and closer to $2.8 T^*$ ``on average''); see
\cite[ \S5.2]{propp-wilson:exact-sampling}.  If we had stepped back by
one each time, then the number of simulated steps would have been
$\binom{T^*}{2} \approx \frac12 (T^*)^2$.  If we had stepped back
quadratically, then the number of simulated steps would have been
about $\frac13 (T^*)^{3/2}$.  For this reason some people prefer
quadratic backoff when $T^*$ is fairly small \citep{moller:personal}.

\q So you need the state space to be monotone?
\ans That would certainly be very useful, but there are many counterexamples to the proposition that is necessary.  See \sref{coupling}.

\q Can you do CFTP on Banach spaces?
\ans I don't know.  Depends on whether or not you can figure out the state at time 0.

\q Where's the proof of efficiency?
\ans
 In the monotone setting, loosely speaking, CFTP is efficient whenever
 the Markov chain mixes rapidly; see \sref{monotone}.  There also
 proofs of efficiency for certain non-monotone settings.

\qa{(?)} But you still need to analyze the mixing time, since otherwise you won't know how long it will take.
\ans
 Wrong.  The principal advantage of using CFTP (or Fill's algorithm)
 is that you \textit{don't\/} need these \textit{a priori\/} mixing
 time bounds in order to run the algorithm, collect perfectly random
 samples, and carry on with the rest of the research project.  The fact
 that these are \textit{perfectly\/} random samples is icing on the cake.
 (But I still think that mixing time bounds are interesting.)

\q But sampling spin glass configurations is NP-hard.
\ans
 The spin-glass Markov chain that you're using is probably very slowly
 mixing in the worst case.  CFTP will not be faster than the mixing
 time of the underlying Markov chain.

\q Is CFTP like a ZPP or Las Vegas algorithm for random sampling?
\ans
 ``Yes'' for Las Vegas, ``sometimes'' for ZPP.  [A Las Vegas algorithm
 uses randomness to compute a deterministic function.  The running
 time is a random variable, but with probability 1 it returns an
 answer, and when it does so, the answer is correct.  A Monte Carlo
 algorithm in contrast does not guarantee a correct answer.  A problem
 is in ZPP if it can be solved by polynomial expected time Las Vegas
 algorithm.]

\q Could you make a hybrid algorithm, which starts out by doing CFTP, but then does something different if CFTP starts to take a long time?
\ans Yes.  This would be like experiment $C_T$.

\q What's the probability that CFTP takes a long time?
\ans
 The tail distribution of the running time decays geometrically.  More
 can be said if the Markov chain is monotone, see \sref{monotone}.

\q What if I don't want to wait for $10^{70}$ years.  Does this make CFTP biased?  Is it really better than forwards coupling?
\ansa{1}
 The answer to the second question is yes.  How long are you willing
 to wait?  With forward coupling, that's how long you wait.  With CFTP
 your average waiting time is probably much smaller.
\ansa{2}
 The answer to the first question depends on what you mean by
 ``biased''.  No-one disputes that the systematic bias is zero.  The
 so-called ``user-impatience bias'' (the effect of an impatient user
 interrupting a simulation and progress) is a second order effect that
 pertains to most random sampling algorithms that most people would
 not call biased.
\ansa{3}
 If you're genuinely concerned about the quality of your random
 samples, you should first spend time picking a good pseudo-random
 number generator.
\ansa{4}
 Anyone concerned about ``user-impatience bias'' should be equally
 concerned about ``user-patience bias'': if an experimenter run
 simulations until some deadline, and then lets the last one finish
 before quitting, then the resulting collection of samples is biased
 to contain more samples that take a long time to generate.  But if
 the experimenter instead aborts the last simulation (unless there are
 no samples so far, in which case (s)he lets it finish), then the
 resulting collection of samples is unbiased.  See
 \citep{glynn-heidelberger:budget}.

\q Can you quantify the user-impatience bias?
\ans
 If you generate $N$ samples before your deadline, and then average
 some function of these samples, the bias is at most $\Pr[N=0]$.  If
 you do something sensible in the event that $N=0$, the bias will be
 even less.  See \citep{glynn-heidelberger:budget}.

\subsection{Historical remarks and further reading}

Monotone-CFTP was developed in 1994, although related ideas had
appeared in the literature prior to that time.
\cite*{asmussen-glynn-thorisson:unknown-markov} and \cite{lovasz-winkler:unknown-markov} (see also \citep{aldous:unknown-markov})
had given algorithms for generating perfectly random samples from the
steady-state distribution of any finite Markov chain; CFTP is
generally more efficient at this task \citep{propp-wilson:unknown-markov-tree}.
\cite{letac:contraction} noticed that if one composed random maps backwards
in time rather than forwards in time, then applying these maps to a given
state typically led to pointwise convergence rather than just
convergence in distribution.  (\cite{diaconis-freedman:rfuncs} survey
this and related work.)  But the random maps were composed backwards in
time forever, and little attention was given to the question of how or
if one could stop the process and obtain a random sample in finite
time; one researcher in this area expressed surprise and disbelief upon
first learning that this was possible \citep{foss:personal}.
Notable exceptions are the \cite{aldous:tree} /
\cite{broder:tree} algorithm for random spanning trees, which reverses
time when building the tree, and an algorithm for the ``dead leaves
model'' in which leaves fall up from the ground rather than down from
the sky (see \citep{jeulin:dead-leaves},
\citep{kendall-thonnes:geometry}, and
\texttt{http://www.warwick.ac.uk/statsdept/Staff/WSK/dead.html}).  In
both these cases a Markov chain is run backwards in time.
Monotone-CFTP and nearly all subsequent versions of CFTP compose their
randomizing operations forwards in time but starting from ever distant
times in the past.  \cite{johnson:testing} independently studied
monotone couplings, but did not couple them from the past.
Monotone-CFTP may be applied to a surprisingly wide variety of
Markov chains (see \sref{monotone}), and after its success, many people
started looking at other coupling methods that could be used with CFTP
(see \sref{coupling}).

For more in-depth explanations of the ideas described so far, the
reader is referred to \citep{propp-wilson:exact-sampling} or the
expository articles \citep{propp:expository} or
\citep{propp-wilson:CFTP-aug}.  Some people prefer the explanation of CFTP in
\citep{fill:interruptible}.  Subsequent to the writing of this primer on CFTP,
the author became aware of two additional expositions on
perfect sampling with Markov chains: \citep{dimakos:guide-exact} and
\citep{thonnes:primer}.

A more recent development is an algorithm related to CFTP but for
which the Markov chain is run forwards in time and never restarted
further back in the past \citep{wilson:rocftp}.

Additional information on perfect sampling is available at
\texttt{http://dimacs. rutgers.edu/\char126dbwilson/exact/}.

\subsection{Overview of common coupling methods}
\label{sec:coupling}

\subsubsection{Monotone coupling}
\label{sec:monotone}

We saw the method of monotone coupling when looking at the toy example
in \sref{cftp}.  Just as a reminder, the state space comes equipped
with a partial order $\preceq$ with a biggest state $\hat{1}$ and a
smallest state $\hat{0}$ such that $\hat0\preceq x\preceq \hat1$ for
all states $x$.  A randomizing operation on a partially ordered set is
monotone if $x\preceq y$ implies $\phi(x,u) \preceq \phi(y,u)$.  To
test if the randomizing operations determine $X_0$, it is only
necessary to apply them to the two starting values $\hat0$ and
$\hat1$.  Examples of monotone couplings include

\begin{itemize}
\item The Fortuin-Kasteleyn random cluster model with $q\geq 1$.  These can be used to generate random configurations from dimer Ising and Potts models \citep{fortuin-kasteleyn:cluster}.
\item Dimer models on bipartite planar graphs.  These include domino tilings \citep{levitov:equivalence,zheng-sachdev:quantum} and lozenge tilings \citep{blote-hilhorst:roughening}.  See also \citep{conway-lagarias:tiling,thurston:conway,propp:tilings}.
\item Ice models, including the 6-vertex model \citep{vanbeijeren:roughening} and the 20-vertex model.
\item The hard core model on bipartite graphs (random independent sets) \citep*{kim-shor-winkler:independent}.
\item The Widom-Rowlinson model with 2 types of particles \citep*{haggstrom-lieshout-moller:exact-spatial}.  This is actually a special case of the hard core model.
\item Linear extensions of a 2D partially ordered set \citep{felsner-wernisch:linear-extension}.
\item Attractive area interaction point process \citep{kendall:area-interaction}.
\item Certain queueing systems \citep{lund-wilson:storage}.
\item The beach model \citep{nelander:beach}.
\item Statistical analysis of DNA using the ``M1-M$m$'' model \citep*{muri-chauveau-cellier:dna}.
\item A mite dispersal model \citep{straatman:thesis}.
\item Slice samplers \citep*{mira-moller-roberts:slice}.
\item The autonormal distribution (\sref{autonormal}).
\end{itemize}

When using Fill's algorithm, a somewhat weaker notion of monotonicity is sufficient.  Rather than a monotone randomizing operation, as is needed for monotone CFTP, it is sufficient to have a monotone pairwise coupling \citep{fill:interruptible}.  Rather than produce a whole random map which is both monotone and marginalizes to a given Markov chain, it is sufficient to specify how any two given states may be updated together in a monotone fashion. \cite{fill-machida:monotonicity} show that there are Markov chains with a pairwise monotone update rule but with no monotone randomizing operation, but so far there haven't been any such examples where someone wanted to sample from the steady state distribution.

When monotone couplings are used either with CFTP or Fill's algorithm, the running time of the algorithm has been rigorously related to the mixing time of the Markov chain (see \cite{propp-wilson:exact-sampling} and \cite{fill:interruptible}).
For CFTP the relevant notion of mixing time $\Tmix$ is the ``total variation threshold time'', which is what most people mean by the phrase ``mixing time''.  For Fill's algorithm, the relevant notion of mixing time is $\Tsep$, the ``separation distance threshold time''.  We won't define these terms here, but the interested reader can read about them in \citep{aldous-fill:book}.  

Let $\ell$ denote the length of the longest chain within the partially ordered state space.  Then the time to coalescence (or how far back in the past that you need to start) $T^*$ has its expected value bounded by $$ \Tmix/e \leq E[T^*] \leq 2 \Tmix (1+\ln \ell),$$
where $e=2.71828\ldots$ \citep{propp-wilson:exact-sampling}.  For some
examples (including lozenge tilings of a hexagonal region) it is
possible to determine both the mixing time in the coupling time to
within constants, and for these cases the coupling time does not
contain an extra $\log$ factor \citep{fill:move-to-front}
\citep{wilson:lozenge+shuffle}.  But it is also possible to construct
examples for which the $\log$ factor does appear in the coupling time
\citep{lund-wilson:storage}.

The expected coupling time for Fill's algorithm can be bounded by
$$ E[T^*] = \Theta(\Tsep).$$
In fact, the cumulative distribution function of the coupling time plus the separation distance as a function of time add up to the constant function $1$ \cite{fill:interruptible}.  For comparison with CFTP, we note that in general $\Tsep \geq \Tmix$, and that for reversible Markov chains $\Tsep = \Theta(\Tmix)$.

Thus for monotone Markov chains, from a running time standpoint, there is little or no reason not to use one of these two perfect sampling algorithms.

\subsubsection{Anti-monotone coupling and Markov random fields}
\label{sec:random-field}

Mathematicians use the term ``Markov random field'' where physicists use the term ``spin system'' where statisticians use the term ``conditionally specified model''.  A Markov random field is a collection of random variables (or spins) defined at the vertices (or sites) of a graph; the edges of the graph contain information about the correlations between the random variables.  If the values of all the spins except one are specified, then the conditional distribution of the remaining spin is a function only of that spin's neighbors.  One of the more frequently used Markov chains on spin systems is the ``single-site heat bath'', also known as ``Gibbs sampling''.  The Markov chain picks a site, either at random or in sequence, and then randomizes the spin at that site by drawing from its conditional distribution when the remaining spins are held fixed.

A spin system is attractive if there is a partial order on the values of the spins, such that increasing the values of the spins only increases the conditional distribution of the spin at a given site.  Most of the examples of monotone coupling given in \sref{monotone} are in fact instances of attractive spin systems.

A spin system is repulsive if there is a partial order on the values of the spins, such that increasing the values of the spins only decreases the conditional distribution of the spin at a given site.  
``Anti-monotone coupling'' was first used by \cite{kendall:area-interaction} and \cite{haggstrom-nelander:antimonotone} to generate perfectly random samples from repulsive spin systems.  Instances of these repulsive systems include
\begin{itemize}
\item The repulsive area interaction point process \citep{kendall:area-interaction}.
\item The hard-core model on non-bipartite graphs \citep{haggstrom-nelander:antimonotone}.
\item The Ising antiferromagnet \citep{haggstrom-nelander:antimonotone}.
\item Fortuin-Kasteleyn random cluster model with $q<1$ \citep{haggstrom-nelander:antimonotone}.
\item The autogamma distribution \citep{moller:conditional}.
\end{itemize}

With anti-monotone coupling with Gibbs sampling, one maintains at each site an upper bound and a lower bound for the spin at that site.  If there is a smallest and a biggest spin value, then these become the initial lower bound and upper bound; see \sref{unbounded} if there is no smallest or largest spin value.  In contrast with attractive spin systems, there is no particular reason for the configuration consisting entirely of the lower bounds or entirely of the upper bounds to have positive probability.  For instance, for the hard-core model, the configuration consisting of all upper bounds will have particles too close to each other (unless there are no edges), and will thus have probability 0.  But this is irrelevant for our purposes: the bounds on the spin values specify a superset of the possible states that the Markov chain may be in.
When doing a Gibbsian update at a given site, when determining the upper bounds at the neighboring sites are used when determining the new lower bound at that site, and vice versa.

The area interaction point process requires further explanation because there are infinitely many sites.  The configurations tend to be sparse because there are only two possible spin values, and (normally) the spins at all but finitely many sites are of the first type.  \cite{kendall:area-interaction} used his method of dominated CFTP (see \sref{kendall}) together with anti-monotone coupling to sample from this distribution.

For the autogamma distribution, the possible spin values are the non-negative reals.  See \sref{compactify} for remarks about upper-bounding the possible spin values.

\cite{haggstrom-nelander:antimonotone} and \cite{huber:techniques}
independently generalized the approach of anti-monotone coupling on
repulsive systems to more generic Markov random fields.  At each
site, one still maintains the set of possible values of the spin at
that site, but that set can no longer be represented by an interval
specified by a lower bound and upper bound.  The sense of possible
spins at each site collectively define some abstract high-dimensional
box which contains the possible states of the Markov chain.  Updating
the set of possible values of the spin at a given site becomes more
complicated than in the anti-monotone setting, and good coupling
methods are essential to making it work.
The reader is referred to the original articles to see how the following examples are done.
\begin{itemize}
\item Random $q$-colorings \citep{haggstrom-nelander:antimonotone,huber:techniques}.
\item The Widom-Rowlinson model with more than two particle types \citep{haggstrom-nelander:antimonotone}.
\end{itemize}

\subsubsection{Techniques for Bayesian inference}
\label{sec:bayesian}

\cite{murdoch-green:continuous} and \cite{green-murdoch:bayesian}
 developed a number of
techniques that are suited for applying CFTP to sampling from the
(continuous) posterior distributions associated with Bayesian
inference problems.  Here the state space is typically $\R^d$.  At
each time step the algorithm maintains a superset of the possible
states of the Markov chain where the superset is represented by a
finite collection of boxes together with a finite set of points.  It
is not possible to do justice to these coupling techniques within a
short primer on CFTP, so the reader is referred to their original
articles.

\subsection{Coupling \epf}
\label{sec:ex-post-facto}

Here we review ``\epf coupling'', a term introduced by Jim
Fill.  Later we explain the role that \epf coupling plays in Fill's
algorithm (\sref{fill}) and ``coupling into and from the past''
(\sref{kendall}).

In ordinary pairwise coupling, a procedure takes as input two states $x$
and $y$ and produces two states $x'$ and $y'$ so that the transition from $x$
to $x'$ and the transition from $y$ to $y'$ both look like they were
produced from a given Markov chain.
\begin{code}
\item      ($x'$,$y'$) := PairWiseCoupling($x$,$y$)
\end{code}
Usually there are additional constraints on useful couplings, such as
a monotonicity constraint or a contraction property.  In \epf
coupling, somebody else generates the $x'$ according a Markov update on
$x$, and your job is to take $x$, $x'$, $y$, and generate a random $y'$ so that
the distribution of $x',y'$ given $x,y$ is governed by the original
coupling.
\begin{code}
\item       $x''$ := MarkovUpdate($x$)
\item       $y''$ := ExPostFactoCoupling($x$,$y$,$x'$)
\item       /* DistributionOf($x''$,$y''$) == Distribution($x'$,$y'$) */
\end{code}
The same idea applies to random mappings.
\begin{code}
\item       $F()$ := RandomMap()
\item       /* $F$ now defines a Markov update from any state $x$ */
\end{code}
Somebody else gives you a Markov chain step $(x,x')$, and your job is to
produce a random map $F$ from a given random map distribution, but
conditioned on $F(x)=x'$.
\begin{code}
\item       $x'$ := MarkovUpdate($x$)
\item       $F'()$ := ExPostFactoMapping($x$,$x'$)
\item       /* DistributionOf($F'$) == DistributionOf($F$) */
\item       /* $F'(x)$ == $x'$                             */
\end{code}
In other words, to do \epf coupling, we generate a random
map or a random pairwise update, conditioned to satisfy a certain
constraint.  Since we just need to sample from a conditional
distribution, in principle any coupling can be done \epf,
but in practice this can be easier said than done.

\subsection{Fill's algorithm}
\label{sec:fill}

Here we briefly describe Fill's algorithm, and in particular the role
that \epf coupling plays in it.  For an explanation of why Fill's
algorithm works, see \citep{fill:interruptible} or
\citep*{fill-machida-murdoch-rosenthal:fill}.

In Fill's algorithm, a single trajectory of the Markov chain is run
forward in time for some number of steps $X_0,\ldots,X_n$.  This trajectory
is treated as a sample path from the time-reversed Markov chain.  Then
a second trajectory (or in subsequent work, many trajectories) of the
time-reversed Markov chain is coupled to it \epf.  In other words,
the time-reversal of the first trajectory together with all the
new time-reversed trajectories have a joint distribution that is
governed by some pre-specified coupling, conditioned upon the
trajectory from $X_n$ being $X_n,X_{n-1},\ldots,X_1,X_0$.  The state $X_n$ is
returned if each of the time-reversed trajectories coalesced to $X_0$.
Otherwise, the current experiment is discarded, and another (independent)
one may be started.
In the monotone setting, Fill's algorithm only requires \epf
coupling for a monotone pairwise coupling, but in more general
settings, random maps are coupled \epf.

Fill's algorithm is interruptible with respect to a deadline specified
in terms of Markov chain steps, so the corresponding user-patience
bias does not affect it.  Code implementing the algorithm may or may
not be interruptible with respect to a deadline specified in terms of
time.

\subsection{Methods for unbounded state spaces}
\label{sec:unbounded}

Suppose that one has a partially ordered state space together with the
monotone (or anti-monotone) randomizing operation.  When doing
monotone or anti-monotone CFTP, having a top state and bottom state
is important, or at the very least, very useful.  What can one do when
there is no top state?  The unboundedness of the state space can be
similarly problematic when the couplings used do not rely on a partial
order.  In this section we describe the three main techniques that
have been used when dealing with unbounded state spaces.

The method in \sref{compactify} extends the state space, the method in
\sref{murdoch} modifies the Markov chain, and the method and
\sref{kendall} uses two coupled Markov chains, one going backwards in
time in the other going forwards.  For those familiar with the term
``uniform ergodicity'', the method in \sref{compactify} requires
uniform ergodicity, the method in \sref{murdoch} produces a uniformly
ergodic Markov chain starting from one that is non-uniformly ergodic,
in the method and \sref{kendall} works with non-uniformly ergodic
Markov chains.  Despite the differences in approach and capabilities
of the methods described in \sref{compactify} and \sref{kendall}, they
are both frequently referred to by the same term, namely ``dominated
CFTP''.  The method in \sref{murdoch} is comparatively new, so it is
too early to tell whether or not it too will be referred to by this
same term.  We mention a fourth method in \sref{multistage}.

\subsubsection{Compactifying the state space}
\label{sec:compactify}

Adjoin a top state or bottom state if these are missing.  Then the state spaces no longer unbounded. Trite as this solution may sound, in more than one case it works just fine and solves the problem, and it is much simpler than the approaches in \sref{kendall} and \sref{murdoch} for dealing with unbounded state spaces.

Let us denote the newly adjoined top and/or bottom states by $+\infty$
and $-\infty$ respectively.  Let $P_x(\cdot)$ denote the probability
distribution of the next state of the Markov chain when it starts in
state $x$.  If there is some probability distribution which
stochastically dominates $P_x(\cdot)$ for each $x$, then in the
monotone case we can set $P_{+\infty}(\cdot)$ to be this distribution,
and in the anti-monotone case we can set $P_{-\infty}(\cdot)$ to be
this distribution.  Similarly, if there is some probability
distribution which is stochastically dominated by $P_x(\cdot)$ for
each $x$, then we can define $P_{-\infty}(\cdot)$ or
$P_{+\infty}(\cdot)$ in the monotone or anti-monotone cases
respectively.  If we adjoined $\pm\infty$ but then were unable to
define $P_{\pm\infty}(\cdot)$, then one of the other methods (in
\sref{murdoch} or \sref{kendall}) for dealing with unbound state
spaces should be used.

In the new Markov chain the states $\pm\infty$ are transient, so the new steady-state distribution is the same as the old one, and we can proceed to sample from it using monotone or anti-monotone CFTP.

This approach is sometimes even easier done than said.  For instance,
when doing anti-monotone coupling with the autogamma distribution, the
system is a repulsive spin system where the possible spin values are
$\R^+$.  Each spin variable $x_i$, conditional upon the remaining spins,
is governed by a gamma distribution with shape parameter $\alpha_i$
which is then scaled down by a factor of 
\begin{equation}
\label{eq:scale}\tag{*}
\beta_i+\sum_{j:j\neq i} \beta_{i,j} x_j,
\end{equation}
where $\beta_i>0$ and $\beta_{i,j}\geq 0$ \citep{moller:conditional}.
On any modern computer we can simply set the top state to be
\texttt{1.0/0}.  This is because all modern computers conform to the
IEEE 754 floating point arithmetic standard, which has built-in
representations for both $+\infty$ and $-\infty$, and knows how to
sensibly add numbers to infinity and divide numbers by infinity; see
\citep{goldberg:arithmetic}.  No special code needs to be written to
sample from $P_{+\infty}$ or otherwise deal with such a large top
state: the code which computes the inverse scale parameter given
by \eqref{eq:scale} and updates the range of possible spin values at a
given site when the neighboring spin values are bounded by finite
values will also work correctly when the neighboring spin values are
bounded by $+\infty$ (thanks to IEEE arithmetic).

(\cite{moller:conditional} did not regard $+\infty$ to be a valid spin
value, and used pages of detailed calculations to
verify the anti-monotone CFTP still works.  When we regard $+\infty$
as a valid spin value, it is obvious without calculation that
anti-monotone CFTP still works.)

In other applications the computer hardware may not come prewired to
deal with the $\pm\infty$ configurations as it did in the autogamma
example, in which case this must be done in software.  When figuring
out whether or not the randomizing operations given by
$U_{-T},U_{-T+1},\ldots,U_{-1}$ determine the state at time $0$, one
piece of code could deal with the random map specified by $U_{-T}$,
and another piece of code could deal with the subsequent $T-1$ random
maps.  From a mathematical standpoint there is no difference between
the first randomizing operation and the subsequent ones.  From an
implementation standpoint, it is sometimes easier to write one piece
of code optimized for the special case of the upper bound being
$+\infty$ (and/or lower bound being $-\infty$), and a separate piece
of code optimized for finite upper and lower bounds.

We remark that the continuous Widom-Rowlinson model is another example
where the state space can be compactified by adjoining a top state
$+\infty$ (consisting of all red points) and a bottom state $-\infty$
(consisting of all blue points).
\cite*{haggstrom-lieshout-moller:exact-spatial} identified a finite
``quasi-maximal'' state \texttt{big} and a finite ``quasi-minimal''
state $-$\texttt{big} such that $P_\infty = P_{\texttt{big}}$ and
$P_{-\infty} = P_{-\texttt{big}}$.  Therefore they were able to
represent $\pm\infty$ using $\pm\texttt{big}$ within their
monotone-CFTP code.  (Currently M\o ller advocates the approach, if not
perspective, of the previous paragraph.)

\subsubsection{Murdoch's method of mixing with an independence sampler}
\label{sec:murdoch}

This method is the next one to try if compactifying the state space
does not work.  This happens when the probability distribution has
infinite tales, and a Markov chain started sufficiently far out in the
tales can take arbitrarily long to reach the main part of the state
space were the steady-state distribution $\pi$ is principally
supported.  The idea is to mix the given Markov chain with one that is
fairly rapid far out in the tales.  \citep{murdoch:exact-bayesian} recommended mixing the
given Markov chain with an ``independence sampler''.  The independence
sampler does a Metropolis-Hastings update, but where the proposal
distribution $P_x()$ starting from state $x$
is independent of $x$.  Normally an
independence sampler by itself will have very poor mixing time
characteristics within the main part of the state space, but this
doesn't matter, since we still use the given Markov chain.  The reason
for using the independence sampler is that when the proposal
distribution has suitably fat tails, {\it all\/} of the states
suitably far out in the tails will in fact get updated.  By ``suitably fat tails'', we mean that the proposal density $P()$ satisfies $P(x)/\pi(x)$ grows as $x\rightarrow\infty$.  If the starting state is $A$ and the proposed state is $B$, then the proposal is accepted with probability
$$\min\left\{1,\frac{\pi(B)P_B(A)}{\pi(A)P_A(B)}\right\} = \min\left\{1,\frac{\pi(B)}{P(B)}\frac{P(A)}{\pi(A)}\right\},$$
which will be $1$ for $A$ suitably far out in the tails of the distribution.
After one
step of the independence sampler, there is some finite box containing
the updated state.  From there we can do coupling with the given
Markov chain.  We give a concrete example of this method in \sref{autonormal} of this article; further examples are given by
\cite{murdoch:exact-bayesian} and \cite{wilson:rocftp}.

\subsubsection{Kendall's method of dominated CFTP / CIAFTP}
\label{sec:kendall}

``Coupling into and from the past'' is an extension of ``coupling from
the past'' introduced by \cite{kendall:area-interaction}, though he used the term ``dominated CFTP'' (see remark below).
In it we have two Markov chains, we already know how to sample from
the stationary distribution of the first chain (the reference chain), and we want to sample
from the stationary distribution of the second chain (the target chain).  It is assumed
that there is a ``useful'' coupling that updates a single state of the
reference chain together with all possible states of the target chain.  A
draw from the stationary distribution of the reference Markov chain is
produced, and then this chain is run backwards into time (via running
the time-reversal forwards in time), producing a sample path of the
reference chain up to time $0$.  Then random maps for the target Markov
chain are randomly generated so that they are coupled \epf to
the sample path of the reference Markov chain.  If we can determine that
there is only one possible value for the state of the target Markov
chain at time $0$, then (assuming we can do this with probability $1$) this
state is a draw from the stationary distribution of the target chain.
Observe that the state of the reference Markov chain at any given time
contains implicit information about the random mappings of the target
Markov chain at all previous times.  This implicit information can be
taken into account when determining the possible states of the target
Markov chain at time 0.  Making use of this implicit information about
previous not-yet-generated random maps is what distinguishes ``coupling
into and from the past'' from ordinary CFTP, and is what enables it to
generate perfectly random samples from non-uniformly ergodic Markov
chains.

We give pseudocode below to make it easier to compare and
contrast ``coupling into and from the past'' with CFTP.

Coupling from the past:
\begin{code}
\item   $T := 1$
\item   \Repeat \{
\item\s   Set := \sspace
\item\s   \For $t$ := $T$ \DownTo $1$
\item\s\s   \If $t$ is a power of 2
\item\s\s\s   SetRandomSeed(seed[$\log_2(t)$])
\item\s\s   ApplyRandomMap(Set)
\item\s   $T := 2*T$
\item   \} \Until Singleton(Set)
\item   \Output ElementContainedIn(Set)
\end{code}

Coupling into and from the past:
\begin{code}
\item  $X[0] :=$ ReferenceChainRandomState()
\item  $T := 1$
\item  \Repeat \{
\item\s  SetRandomSeed(seed1[$\log_2(T)$])
\item\s  \For $t := \lfloor T/2 \rfloor+1$ \To $T$
\item\s\s  $X[t]$ := ReverseReferenceChain($X[t-1]$)
\item\s  Set := $\langle$portion of state space compatible with $X[T]\rangle$
\item\s  \For $t := T$ \DownTo $1$
\item\s\s  \If $t$ is a power of $2$
\item\s\s\s  SetRandomSeed(seed2[$\log_2(t)$])
\item\s\s  ApplyTargetChainRandomMapCoupledExPostFacto($X[t]$,$X[t-1]$,Set)
\item\s  $T := 2*T$
\item  \} \Until Singleton(Set)
\item  \Output ElementContainedIn(Set)
\end{code}
     
The above description may seem abstract, but
\cite{kendall:area-interaction} gives a concrete example carrying out
all these ideas.  The long awaited article by
\cite{kendall-moller:exact-spatial}, and the article by
\cite{lund-wilson:storage}, give more examples of coupling into and
from the past.
Later we will return to the algorithm in \cite{lund-wilson:storage},
since the (\epf) coupling methods described in \sref{exponential} and
\sref{lmc-epf} significantly simplify it.

\noindent\textit{Remark:}
\cite{kendall:area-interaction}
originally referred to his method by ``dominated CFTP''
because of the role that stochastic domination place in the examples
that he gave.  We prefer the term ``dominated CIAFTP'' or ``CIAFTP''
for two reasons: (1)
``dominated CFTP'' is ambiguous since it by now also refers to the
method in \sref{compactify}, and (2) there is the least one instance
of CIAFTP for which there is no partial order or stochastic domination
\citep{wilson:rocftp}; ``dominated CFTP'' would be a misnomer for this
case, and ``CIAFTP'' sounds better than ``undominated dominated CFTP''.

\subsubsection{A multistage method}
\label{sec:multistage}

 \cite{murdoch:exact-bayesian} proposed that a multistage version of
CFTP due to \cite{meng:multistage-backwards} could be adapted to
sample from unbounded state spaces.  For the application that Murdoch
considered, he found that his other method of mixing with an
independence sampler (\sref{murdoch}) worked better.  The interested
reader is referred to \cite{murdoch:exact-bayesian} for further
information.

\section{Multishift Coupling}

\subsection{Introduction}

A multishift coupler generates a random function $f(x)$ so that for
each real number $x$, the random number $f(x)-x$ is governed by the
same fixed probability distribution, independent of $x$.  A multiscale
coupler is defined similarly, except that $f(x)/x$ is governed by the
same distribution for each positive $x$.  The trivial multishift
coupler, say for the normal distribution, would pick a normally
distributed random variable $X$, and set $f(x) = x+X$ for each real
$x$.  An obvious property of this coupling is that regardless of $X$,
each real number is in the image of $f()$, i.e.\ $f({\mathbb
R})={\mathbb R}$.  \cite{green-murdoch:bayesian} devised a more
sophisticated multishift coupler, the ``bisection coupler'', whose
image is a {\it discrete\/} set of points.  In Green and Murdoch's
application, where a computation is done for each point in the image
under $f()$ of a finite interval, the discreteness of the image is
vital.  But while the number of points in the image is finite for the
bisection coupler, the expected number is infinite.  We develop here
the class of layered multishift couplers, which have more pleasant
properties.  For the standard normal distribution, for instance, our
multishift coupler maps an interval of length $\ell$ to fewer than
$2+\ell/2.35$ points.  Our multishift couplers are also monotone,
i.e.\ $f(x_1)\leq f(x_2)$ when $x_1\leq x_2$, a property not enjoyed
by the bisection coupler.  Monotonicity has proved to be very useful
in a multitude of recent sampling algorithms.  In addition to making
Green and Murdoch's application easier, using these monotone
multishift and multiscale couplers, we develop in \sref{autonormal} an
algorithm for generating perfectly random samples from the autonormal
distribution, improve an algorithm of
\cite{moller:conditional} for sampling from the autogamma
distribution, and simplify the algorithm of \cite{lund-wilson:storage}
for sampling from the stationary distribution of certain storage
systems.

All of these applications involve algorithms based on coupling from
the past.  In each case the Markov chain draws a point from a
distribution which is shifted by a different amount depending on the
starting state, so in one way or another some form of multishift
coupling is used.  When running CFTP it is desireable to use a
randomizing operation that maps large numbers of states to the same or
nearby values --- which should explain in part why it is desirable for
a multishift coupling to have a discrete image.

\subsection{Comparison of multishift couplers}

\ \\
\begin{tabular}{|l||c|c|c|c|}
\cline{2-5}
\multicolumn{1}{c}{} & \multicolumn{4}{|c|}{multishift coupler}\\
\cline{2-5}
\multicolumn{1}{c|}{} & trivial   & Poisson  & bisection & layered \\
\cline{2-5}
\hline
\begin{tabular}{@{}l@{}}
works for which\\
distributions?
\end{tabular} 
                    & all    & exponential & \begin{tabular}{@{}c@{}}symmetric\\ unimodal\end{tabular} & \begin{tabular}{@{}c@{}}nonsingular\\univariate\end{tabular} \\
\hline
discrete image?    & no       & yes       & yes       & yes     \\
\hline
\begin{tabular}{@{}l@{}}
expected size of\\
image of finite region
\end{tabular} & uncountable&finite   & typically $\infty$ & \begin{tabular}{@{}c@{}}
typically finite\\
(see \sref{expected-size})
\end{tabular}\\
\hline
\begin{tabular}{@{}l@{}}
number of parameters\\
specifying coupling
\end{tabular}  & 1      & $\infty$  & 2         & 3       \\
\hline
monotone?            & yes    & yes       & no        & yes     \\
\hline
has been used for    & autogamma & dams    & posteriors & autonormal \\
\hline
see remarks in       & \sref{autogamma} & \sref{dams} & \sref{posteriors} & \sref{autonormal-} \\
\hline
\end{tabular}

\subsection{Applications of multishift coupling}

\subsubsection{Autogamma (pump reliability)}
\label{sec:autogamma}

\cite{moller:conditional} proposed a CFTP-based algorithm for
sampling from the ``autogamma'' distribution (defined in \sref{compactify}), which
governs the posterior distribution of the pump reliability problem of
\cite{gelfand-smith:sampling}.  Previously
\cite{murdoch-green:continuous} had applied their techniques to obtain
a CFTP-based algorithm for this problem; M\o ller's approach was more
specialized and efficient.  The
output produced is numerically within a user-specified $\varepsilon$
from an ideal exact output that has zero bias.  In his paper, and also
at two recent conferences, M\o ller pointed out that some sort of
hybrid algorithm, the algorithm he described joined with one of the
\cite{murdoch-green:continuous} methods, could reduce the numerical
error $\varepsilon$ to zero.  \cite{moller:conditional} also described
another way that $\varepsilon$ could be reduced to 0, but noted that
the method was not practical.  Using the multishift coupler described
in \sref{gamma} for the gamma distribution, only a few small changes to M\o
ller's algorithm are needed to drive the error $\varepsilon$ to zero.
When M\o ller's code is so modified, not only do we acheive the
theoretically pleasing $\varepsilon=0$, but the running time is
slashed as well.
 \cite{moller:conditional} reported the following empirical expected ``$\varepsilon$-coalescence''
 times associated with various values of $\varepsilon$:
\begin{tabular}{|l||c|c|c|c|c|c|}
\hline
accuracy $\varepsilon$              & $10^{-3}$ & $10^{-4}$ & $10^{-5}$ & $10^{-8}$ & $10^{-14}$ & ``$0$''; \begin{tabular}{@{}l@{}}machine\\precision\end{tabular} \\
\hline
\begin{tabular}{@{}l@{}}expected time\\ for $\varepsilon$-coalescence\end{tabular} & $9.3047$  & $11.3170$ & $13.3262$ & $19.3508$ & $31.3775$  & $34.8263$ \\
\hline
\end{tabular}

\noindent
Using the layered multishift coupler for the gamma distribution (\sref{gamma}) we get $\varepsilon=0$ and an empirical expected coalescence time of $5.219$.

\subsubsection{Storage systems}
\label{sec:dams}
The algorithm of \cite{lund-wilson:storage} for sampling from the
steady-state distribution of certain storage systems used a multishift
coupler for the exponential distribution, but the random function
$f()$ output by the coupler required infinitely many parameters to
specify it.  Nonetheless only finitely many of these parameters were
required to evaluate the function $f()$ at finitely many points, so by
generating the requisite parameters on the fly, the computation could
be kept finite.  But generating these parameters in a consistent and
time- and space-efficient manner, while still allowing the CFTP
protocol to re-read the same random map when it needs to, was not
entirely trivial.  The layered multishift coupler for the exponential
distribution (\sref{exponential}) (indeed, for any distribution) only
requires three parameters to specify the random function $f()$.  Using
this coupling offers a significant conceptual and coding
simplification.

\subsubsection{Bayesian inference techniques}
\label{sec:posteriors}
For purposes of doing Bayesian inference,
\cite{green-murdoch:bayesian} generate random samples using a CFTP
algorithm based on the Metropolis-Hastings update rule.  If the
current state of the system is $(x_1,\ldots,x_n)$, a proposed new
state $(y_1,\ldots,y_n)$ is generated according to a normal
distribution centered about $(x_1,\ldots,x_n)$.  Then an appropriately
weighted coin is flipped to determine whether the next state should be the
old state $(x_1,\ldots,x_n)$, or the proposed state
$(y_1,\ldots,y_n)$.  Since they use their bisection coupler for the
normal distribution, they can perform this Metropolis-Hastings update
rule starting from all states (in a finite portion of $\R^n$), and the
set of proposed states will be a finite set.  But as mentioned before,
with the bisection coupler the set of proposed points can be very
large, and its expected size is infinite.  \cite{green-murdoch:bayesian} dealt with this feature without producing an
algorithm with infinite expected running time.  Using instead the layered
multishift coupler would simplify the
algorithm, since there no longer needs to be code to deal with the possibility of a very large image, and could 
make the algorithm more efficient for higher dimensional problems,
as explained in \sref{normal-d}.

\subsubsection{Autonormal}
\label{sec:autonormal-}

In \sref{autonormal} we see how to apply CFTP to perfectly sample from
the autonormal distribution.  There we use monotone-CFTP, so it is
important for the multishift coupler to be monotone, which rules out
the bisection coupler.  If sampling from the autogamma is any
indication of what would happen with the autonormal, using the trivial
multishift coupler would be both unexact and inefficient.  For this
reason we use the layered multishift coupler of \sref{normal}
in \sref{autonormal}.

\subsection{Layered Multishift Coupler}
\label{sec:multishift}

\subsubsection{Rectangular distribution}

We warm up with the rectangular distribution.  The multishift coupler
for the rectangular distribution is illustrated in \fref{rectangle},
and is given algebraically below.

\begin{figure}[phtb]
\centerline{\epsfig{figure=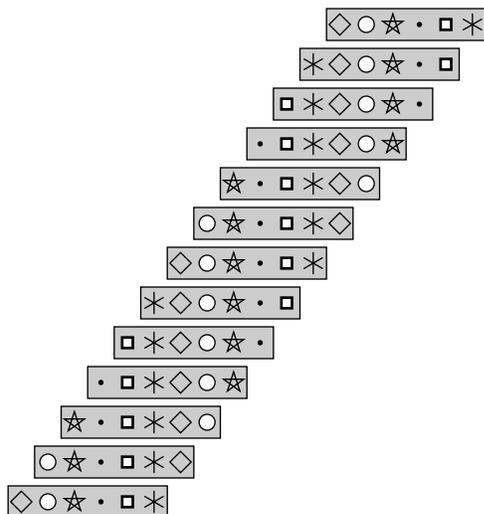}}
\caption{
Illustration of multishift coupling for the rectangular distribution.
A random point $X$ is drawn from the rectangular distribution with
endpoints $L$ and $R$; in the figure six possible such points are
denoted by six different symbols.  If the square, for instance, is
chosen in the first rectangle, then for any other shifted version of
the rectangular distribution, the square will be chosen.  If a
uniformly random symbol/point is drawn from one of the rectangles,
then the corresponding symbol/point in any other given fixed rectangle
will also be uniformly random.
}
\label{fig:rectangle}
\end{figure}

\begin{coupling}
\item Parameters:
\item \s  $L =$  left endpoint of rectangle
\item \s  $R =$ right endpoint of rectangle
\item Random variables:
\item \s   $X := \Uniform(L,R)$
\item Mapping:
\item \s   \fslrx
\end{coupling}

Since $(R-X)/(R-L)$ is uniformly distributed between $0$ and $1$, for any
fixed $s$, the fractional part of $(s+R-X)/(R-L)$ is uniformly
distributed between $0$ and $1$.  If we ignored the floor in the
definition of $f(s)$, the expression would simplify to $s+R$.  If we refrain from ignoring the floor,
then a uniformly random quantity between $0$ and $(R-L)$ is subtracted
from $s+R$, so that $f(s)$ is uniformly distributed between $s+L$ and
$s+R$, as desired.

It is the floor that makes the image of $f_{L,R,X}$ discrete.  If we
continuously increase $s$ by $R-L$, the value of $f_{L,R,X}(s)$
changes only once.  It is also clear that $f_{L,R,X}(s)$ is monotone in $s$.

\subsubsection{Normal distribution}
\label{sec:normal}

Since we already know how to do multishift coupling for the
rectangular distribution, to do the normal distribution, we just need
to express it as a convex combination of rectangular distributions.
If we pick a random rectangle according to a suitable distribution,
and then pick a random point within the rectangle, then the result is
a normally distributed random variable.  What we will do is pick a
random rectangle according to the suitable distribution, and then do
multishift coupling with the corresponding rectangular distribution.

There is more than one way to express the normal distribution as a
convex combination of rectangular distributions, and we illustrate two
of these in \fref{layer}.  In the left part of the figure we have stated the
region between the $x$-axis and the probability density function of
the normal distribution.  (We choose not to normalize the density
function by $1/(\sqrt{2\pi}\sigma)$, since the method works fine with
unnormalized density functions.)  It is well known that if we pick a
uniformly random point $(X,Y)$ from the region, then $X$ will be
distributed according to a normal distribution.  In the right part of
the picture we have taken the portion of the region lying to the left
of the $y$-axis and reflected it vertically about the line $y=1/2$.
If the pick a uniformly random point $(X,Y)$ from this modified
region, then $X$ will still be distributed according to a normal
distribution.  We may view each of these regions as being composed of
a stack of many very thin horizontal rectangles (or ``layers''), as
shown in the lower panels of the figure.

\begin{figure}[phtb]
\centerline{\epsfig{figure=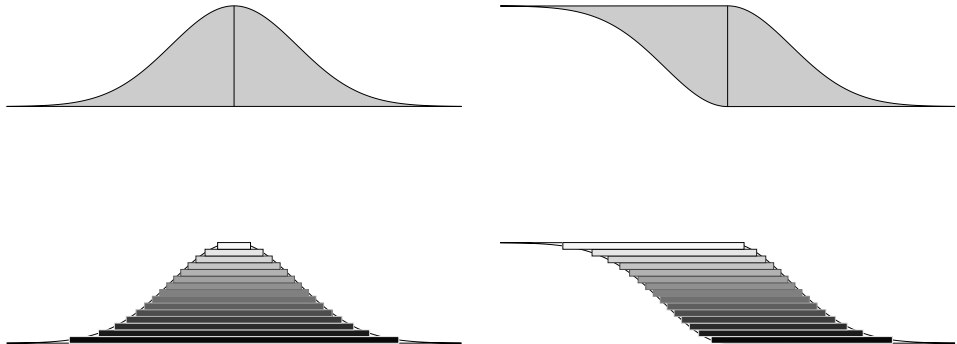,width=\textwidth}}
\caption{
Two different ways of expressing the normal distribution as a convex
combination of rectangular distributions.  A given rectangle is chosen
with probability proportional to its length.  The region on the right
is obtained by vertically reflecting the left portion of the region on
the left.}
\label{fig:layer}
\end{figure}

Let $L$ and $R$ denote the ($x$-coordinates of the) left and right
endpoints of the rectangle containing the uniformly random point
$(X,Y)$ in the region.  If we condition the point $(X,Y)$ to lie
within a particular rectangle, then its distribution within the
rectangle is uniformly random.  In particular, $X$ is uniformly random
between $L$ and $R$.  If we let $L$ and $R$ specify the random
rectangle, and we let $X$ be the uniformly random point from the
corresponding rectangular distribution, then since we already know
that $X$ is distributed according to a normal distribution, we have
our desired decomposition of the normal into a convex combination of
rectangular distributions.

We choose to decompose the normal distribution into rectangular
distributions using the region on the right of \fref{layer} rather than
the region on the left, because for the region on the left there is no
lower bound on how short the rectangles can get, whereas for the
region on the right, there is a positive minimum length for the
rectangle.

To pick the point $(X,Y)$ from the region, we draw $X$ from the normal
distribution, determine the range of possible values for $Y$, and then
pick $Y$ uniformly at random from within this range.  To compute $L$
and $R$, we need to be able to invert the probability density
function, which we can easily do for the normal distribution.  This
procedure is re-expressed below.

\begin{coupling}
\item Parameters:
\item \s   $\sigma =$ standard deviation of normals
\item Random variables:
\item \s   $X := \sigma\Normal(0,1)$
\item \s   $Y := \exp(-(X/\sigma)^2/2) \Uniform(0,1)$
\item \s   If ($X<0$) Then $Y := 1-Y$
\item \s   $L := -\sigma\sqrt{-2\log(1-Y)}$
\item \s   $R :=  \sigma\sqrt{-2\log(Y)}$
\item Mapping:
\item \s   \fslrx
\end{coupling}

Regarding the efficiency of the coupling, it is clear that either $L$ or
$R$ will be at least $\sigma\sqrt{\log 4}$ in absolute value.  To show that
$R-L \ge \sigma 2\sqrt{\log 4}$, we note
  $$w(y) = \sqrt{-2\log y} + \sqrt{-2\log(1-y)}$$
is analytic on $(0,1)$, diverges to infinity as $y\rightarrow 0$ and $y\rightarrow 1$, and show that
$w'(1/2+t)$ has a real zero only at $t=0$.  It follows that $w(1/2)$ minimizes $w$
on $(0,1)$.
\begin{align*}
          w'(1/2+t) = {1/2\over \sqrt{-2\log(1/2+t)}} {-2\over 1/2+t}
                   +& {1/2\over \sqrt{-2\log(1/2-t)}} { 2\over 1/2-t} = 0\\
  -2\log(1/2+t)(1/2+t)^2 =& -2\log(1/2-t)(1/2-t)^2
\end{align*}

Since $\log(1/2+t)$ and $(1/2+t)^2$ are strictly monotone increasing on
$(-1/2,1/2)$, the above equation can have at most one root, which we
know must be at $t=0$.  Thus the minimum width (normalized to $\sigma$)
of any rectangle used in the above procedure is $w(1/2) = 2\sqrt{\log 4}$,
which is about $2.35482$.

From this it follows that an interval of length $\ell$ is mapped under
$f_{L,R,X}$ to at most $\lceil 1+\ell/(2.35\sigma) \rceil$ points.

\subsubsection{Multidimensional normal distribution}
\label{sec:normal-d}

Extending this approach to the spherically symmetrical
multidimensional Gaussian distribution is easy, since the coordinates
are independent and are individually Gaussian.  Space is divided into
rectangular regions, each of which is mapped to a single point.  The
volume of each rectangular region is at least $(2.35
\sigma/\sqrt{d})^d$.  Thus when the Gaussian is used as a proposal
distribution for a Metropolis-Hastings update, and storage is required
for each point in the image, even a modest improvement in the constant
factor can be significant for large $d$.

\subsubsection{Unimodal distributions}
\label{sec:unimodal}

The astute reader will have noticed that the only properties about the
normal distribution that we used in the coupling procedure of
\sref{normal} are that we can sample from it, it is unimodal (and we
know where the mode is), and the probability distribution function
$\PDF()$ and its inverse are easy to compute.  For other distributions
with these properties we can use essentially the same procedure for
multishift coupling:

\begin{coupling}
\item Random variables:
\item\s  $X := $ RandomSampleFromDistribution()
\item\s  $Y := \PDF(X) * \Uniform(0,1)$
\item\s  If ($X<\Mode$) Then $Y := \PDF(\Mode)-Y$
\item\s  $L := \LIPDF(\PDF(\Mode)-Y)$
\item\s  $R := \RIPDF(Y)$
\item Mapping:
\item \s   \fslrx
\end{coupling}

If the probability distribution has a density function $\PDF()$ with a
single mode, then the inverse PDF to the left of the mode is
well-defined\footnote{Well-defined almost everywhere.  One might worry
about PDF's with lots of horizontal portions, but the selected $Y$ is
almost always a value where the inverse PDF is well-defined.}, and
similarly to the right of the mode.  This does not necessarily mean
that we can compute these inverses effectively (cf.\ \sref{gamma}), but
in a certain abstract sense it means that layered multishift can be
applied to any unimodal distribution.

Unless the entire distribution is to one side of the mode (which
happens in \sref{exponential}), there will be a positive minimum
length for the rectangles.  Thus not only will the image of a finite
interval have finite expected size, but the image size will be
deterministically bounded.

The alternative generalization (given below) to unimodal distributions
is also noteworthy, in that it provides a ``maximal coupling'': for
any two $s_1$ and $s_2$, the probability that $f_{L,R,X}(s_1) =
f_{L,R,X}(s_2)$ is at least as large as it would be for any other
coupling of the two distributions.

\begin{coupling}
\item Random variables:
\item\s  $X := $ RandomSampleFromDistribution()
\item\s  $Y := \PDF(X) * \Uniform(0,1)$
\item\s  $L := \LIPDF(Y)$
\item\s  $R := \RIPDF(Y)$
\item Mapping:
\item \s   \fslrx
\end{coupling}

\subsubsection{Multimodal distributions}
\label{sec:multimodal}

Even more generally, suppose that the probability density function has
multiple modes, and let us assume that the PDF is well-behaved (e.g.\
almost everwhere differentiable).  We don't give pseudocode for this
case, but it is easy to describe in words.  Refer back to
\fref{layer}, and recall that we did a vertical reflection to one side
of the mode.  For multimodal distributions, one could simply reflect
the region at each place where the derivitive changes sign, and
proceed as before.

\subsubsection{Expected image size}
\label{sec:expected-size}

We start by computing the expected image size of an interval when the
layered multishift coupler is applied to a unimodal distribution.  The
same formula will hold whether or not we vertically reflect the region
on the PDF at its mode.  Suppose that a rectangle with endpoints $L$
and $R$ is selected, let $W=R-L$ denote its width.  Then conditional
on this rectangle being selected, the expected image size of an
interval of length $\ell$ is $1+\ell/W$.  Thus the (unconditional)
expected image size is $1+\ell E[1/W]$.  Let $y$ be the vertical
coordinate of a thin rectangle with length $W$.  The probability that
this rectangle is selected is $W dy$.  Thus $$E[1/W] =
\int_{y=0}^{y=y_{\max}} [1/W]W dy = y_{\max}$$ where
$y_{\max}=\PDF(\Mode)$ is the height of the distribution at its mode.
For instance, using either of our multishift couplers for the normal
distribution, an interval of length $\ell$ is mapped under $f_{L,R,X}$
to on average $1+\ell/[\sqrt{2\pi}\sigma] \doteq
1+\ell/(2.5066\sigma)$ points.

In the case of multimodal distributions, if we reflect the region
under the PDF each time the derivitive changes sign, then the same
reasoning used above still works, except that now
$$y_{\max} = \sum_{\text{local maxima $x$}} \PDF(x) - \sum_{\text{local minima $x$}} \PDF(x).$$

\noindent\textit{Remark:}
If instead of measuring expected image size of $f()$, we measured the
expected number of times that $f(x)$ changes as $x$ increases, then
for unimodal distributions the layered multishift coupler is optimal
in that it minimizes the expected number of changes in $f(x)$.
Consequently the layered multishift coupler (for unimodal
distributions) also has smallest expected image size among the class
of monotone multishift couplers.

\subsubsection{Exponential distribution}
\label{sec:exponential}

Macro-expanding the generic unimodal procedure we get

\begin{coupling}
\item Parameters:
\item \s   $\mu =$ mean of exponential
\item Random variables:
\item \s   $X := \mu \Exponential(1)$
\item \s   $Y := \exp(-X/\mu) \Uniform(0,1)$
\item \s   $L := 0$
\item \s   $R := \mu (-\log(Y))$
\item Mapping:
\item \s   \fslrx
\end{coupling}

which we can simplify to

\begin{coupling}
\item Parameters:
\item \s   $\mu =$ mean of exponential
\item Random variables:
\item \s   $X_1 := \mu \Exponential(1)$
\item \s   $X_2 := \mu \Exponential(1)$
\item Mapping:
\item \s   $\displaystyle f_{X_1,X_2}(s) = \left\lfloor \frac{s+X_2}{X_1+X_2}\right\rfloor (X_1+X_2) +X_1$
\end{coupling}

In \sref{gamma} we will use the observation that if we subtract $X_2$
rather than add $X_1$, then $f(s)$ will be distributed as $s -$
(rather than $+$) an exponential with mean $\mu$.

Since the entire exponential distribution is to the right of its mode,
we no longer have a deterministic upper bound on the size of the image
of a finite interval.  But from \sref{expected-size} we see that the
image of an interval of length $\ell$ will have expected size
$1+\ell/\mu<\infty$.  We remark that this expected image size is equal
to that of the Poisson multishift coupler used by
\cite{lund-wilson:storage}.

\subsubsection{Scaled gamma distribution}
\label{sec:gamma}

Rather than ask that $f(s)$ be distributed as $s +
\langle\text{reference distribution}\rangle$, one could instead ask for the
distribution to be $s \times \langle\text{reference
distribution}\rangle$.  This multiscale coupling can of course be
reduced to multishift coupling of $\log(f(s))$, so our above
techniques can be applied.

One distribution that has been multiscaled in this way \citep{moller:conditional} is the gamma
distribution, which includes as a special case the exponential distribution.
Recall that a gamma random variable with shape parameter $\alpha$ and
scale parameter 1 has a probability density function given by
 $$\PDF(x) = x^{\alpha-1} e^{-x} / \Gamma(\alpha).$$
As mentioned earlier, since \cite{moller:conditional} used the coupling $f_U(s) = s U$ where $U$ is a gamma
random variable, the image of $f_U()$ is the continuum.

We could be methodical and specialize the layered multishift coupler
for unimodal distributions.  Inverting the PDF would require us to
solve a transcendental equation, which we would presumably do via
Newton's method.  But there is more than one way to decompose a
distribution into rectangles.  We describe a second method which only
uses the standard elementary functions.  It is this second method that
was used in the timing experiment reported in \sref{autogamma}.

It is well known (in some circles) that if $G$ is a gamma random
variable with shape parameter $\alpha+1$, and $T$ is an independent random
variable with exponential distribution and mean $1/\alpha$, then the
distribution of $G e^{-T}$ is a gamma distribution with shape parameter $\alpha$.
(The reader unfamiliar with this fact can easily verify it by doing
some calculus.)  So if we scale $G$ by $e^{\log(s)-T}$ we will get a gamma
with the desired shape and scale parameters.  Using our above shift coupler
for the exponential distribution with negative mean, we get the
following procedure

\begin{coupling}
\item Parameters:
\item \s   $\alpha =$ shape parameter of gamma distribution
\item Random variables:
\item \s   $G := \GammaD(\alpha+1,1)$
\item \s   $X_1 := \Exponential(1)/\alpha$
\item \s   $X_2 := \Exponential(1)/\alpha$
\item Mapping:
\item \s   $\displaystyle f_{G,X_1,X_2}(s) = G \exp\left[\left\lfloor \frac{\log(s)+X_2}{X_1+X_2}\right\rfloor (X_1+X_2) -X_2\right]$
\end{coupling}

Once this coupling is written down, it is fairly effortless to use it
within a program.  Since we are using our earlier shift coupler for
the exponential distribution, the number of points in the image of a
finite interval will be finite, unbounded, but with finite
expectation.

\noindent\textit{Remark:}
Since the coupling relies on the multishift coupler for the
exponential, one sees that the expected number of points in the image
of an interval with aspect ratio $r$ will be $1+\alpha \log r$.  If we
had instead been methodical and specialized our coupler for unimodal
distributions, a few calculations reveal that the expected image size
would be $1+[\alpha^\alpha e^{-\alpha}/\Gamma(\alpha)] \log r$, or
about $1+\sqrt{\alpha/(2\pi)} \log r$ for large $\alpha$.

\subsection{Layered multishift coupling \epf}
\label{sec:lmc-epf}

One of the hardest parts of using Fill's algorithm is doing the \epf
coupling \citep{murdoch:fields}.  (The reader should read
\sref{ex-post-facto} if (s)he has not done so already.)  \Epf coupling
is also required when doing ``coupling into and from the past''
(\sref{kendall}).  So as to facilitate the use of these algorithms
when multishift coupling is needed, here we see how to do multishift
coupling \epf.  In fact, the algorithm given by
\cite{lund-wilson:storage} for sampling from the water-level
distribution of the infinite dam uses CIAFTP and a multishift coupler
for the exponential distribution.  As a consequence, in order to
substitute the layered multishift coupler, we need to be able to do
the layered multishift coupling \epf.

Somebody else picks some $s_0$, and generates a random variable $X_0$
from the given distribution shifted by $s_0$.  Our job is to
generate a random $f()$ such that
\begin{itemize}
\item $f(s_0) = X_0$
\item When we randomize over the choices of $X_0$,
       the distribution of $f()$ is what it would be if
       we had simply generated it using the methods in \sref{multishift}.
\end{itemize}

\subsubsection{Unimodal distributions}

The appropriate modification of the coupler in \sref{unimodal} for
unimodal distributions is given below.

\begin{coupling}
\item Somebody else does:
\item\s  $X_0 := s_0 +$ RandomSampleFromDistribution()
\item 
\item Random variables we generate:
\item\s  $X := X_0 - s_0$
\item\s  $Y := \PDF(X) * \Uniform(0,1)$
\item\s  If ($X<\Mode$) Then $Y := \PDF(\Mode)-Y$
\item\s  $L := \LIPDF(\PDF(\Mode)-Y)$
\item\s  $R := \RIPDF(Y)$
\item Mapping:
\item \s   \fslrxepf
\end{coupling}

First note that if $s_0=0$ this above modification works: When we
randomize over the choices of $X_0$ that someone else makes, we just
get our previous multishift coupler.  Furthermore, since $(R-X)/(R-L)$
is between $0$ and $1$ (and $X=X_0$), when we evaluate $f_{L,R,X_0}()$ at
$s_0=0$ we get $X_0$, as desired.

If $s_0\neq 0$, then we can define $g(s)=f(s_0+s)-s_0$.  When we
randomize over $X_0$, the statistical properties of $g$ are identical
to those of $f$.  The condition $f(s_0)=X_0$ translates to $g(0)=X$,
so we can do the \epf coupling with $g$.  Then we translate
back in terms of $f$ by $f(s)=g(s-s_0)+s_0$, which simplifies to the
above stated formula.

\subsubsection{Scaled gamma distribution}

In \sref{gamma} we gave an \textit{ad hoc\/} layered multiscale
coupler for the gamma distribution, which had the virtue of not
requiring the ability to compute the inverse probability distribution
function.  For completeness we describe here how to do this coupling \epf.

\begin{coupling}
\item Somebody else does:
\item\s  $G^* := s_0 \times \GammaD(\alpha,1)$
\item 
\item Random variables we generate:
\item\s  $G_{\alpha} := G^* / s_0$
\item\s  $X := \Exponential(1)$
\item\s  $G_{\alpha+1} := G_{\alpha} + X$
\item\s  $X_2 := \log(1+X/G_{\alpha})$
\item\s  $X_1 := \Exponential(1)/\alpha$
\item Mapping:
\item\s  $\displaystyle \begin{aligned}f_{s_0,G^*,X_1,X_2}(s) &= G^* \exp\left[\left\lfloor \frac{\log(s/s_0)+X_2}{X_1+X_2}\right\rfloor (X_1+X_2)\right] \\
          &= G_{\alpha+1} s_0 \exp\left[\left\lfloor \frac{\log(s/s_0)+X_2}{X_1+X_2}\right\rfloor (X_1+X_2) -X_2\right]\end{aligned}$
\end{coupling}

The key observation is that $G_{\alpha+1}$ and $X_2$ are independent
of one another, and that $G_{\alpha+1}$ is gamma variate with shape
parameter $\alpha+1$ and $X_2$ is an exponential variate with mean
$1/\alpha$.  This we leave as a (perhaps nontrivial) exercise to the
reader.  Once this observation is verified, the rest should by now be routine.

\subsection{Possible extensions}

Duncan \cite{murdoch:personal} has suggested an extension of the
layered multishift coupler, which instead of coupling together normals
with different means, couples together normals with both different
means and different variances.

It is natural to investigate how one might couple together other
 multiparameter families of distributions.  For
instance, if one wanted to do simulations of what physicists would
call a ``$\phi^4$ theory'', then rather than couple together normally
distributed random variables with different means, one would want to
couple together random variables whose unnormalized densities are of
the form $\exp(-x^4 +a x^3 + b x^2 + c x)$, and couple these for the
various values of $a$, $b$, and $c$.  In other words, we'd like to
create a random function $f(a,b,c)$ so that the image of $f()$ is
discrete, but such that for each fixed $a,b,c$, the random value
$f(a,b,c)$ has the appropriate distribution.

\pagebreak
\section{Perfect Sampling of Autonormal Distributions}
\label{sec:autonormal}

\subsection{Background}

\subsubsection{Applications in statistics and physics}
\label{sec:autonormal-apps}

The autonormal is an important distribution that arises in both
statistics and physics.
We quote from lecture notes written by Julian Besag:
\begin{quote}
The conditional autoregressive or auto--Normal formulation was proposed in 
Besag (1974, 1975), though it stems from the stationary infinite lattice 
autoregressions of L\'{e}vy (1948) and Rosanov (1967). Gaussian autoregressions
have been used in a wide range of applications, including human geography 
(e.g.\ Cliff and Ord, 1975, 1981, Ch.\ 4), agricultural field experiments 
(e.g.\ Bartlett, 1978; Kempton and Howes, 1981; Martin, 1990; Cressie and 
Hartfield, 1993), geographical epidemiology (e.g.\ Clayton and Kaldor, 1987; 
Marshall, 1991; Molli\'{e} and Richardson, 1991; Bernardinelli and Montomoli, 
1992; Cressie, 1993, Ch.\ 7), astronomy (e.g.\ Molina and Ripley, 1989; 
Ripley, 1991), texture analysis (e.g.\ Chellappa and Kashyap, 1985; 
Cohen et al., 1991; Cohen and Patel, 1991), and other forms of imaging 
(e.g.\ Chellappa, 1985; Jinchi and Chellappa, 1986; Cohen and Cooper, 1987; 
Simonchy et al., 1989; Zerubia and Chellappa, 1989). Generalizations to 
multivariate $X_i$'s are considered by Kittler and F\"{o}glein (1984) and 
by Mardia (1988), in the context of remote sensing.
\end{quote}

In physics the autonormal distribution is called a ``free field'', or
more precisely, a ``discrete free field''.  Certain statistical
mechanical models (such as the 2D Ising model) have limiting
behaviors, in the limit of large system sizes, that are described by
free fields.  See \citep{spencer:free-field} for background on free
fields in physics.

\subsubsection{Definition}

A (discrete) free field is a (autonormal) distribution on $n$ random
``height'' variables $x_1,\ldots,x_n$, with interaction strength
$F_{i,j}\geq 0$ between variables $x_i$ and $x_j$.  (In general some
of the interaction strengths may be negative, and under suitable
conditions the distribution will still be well-defined.  We assume in
\sref{using-murdoch} non-negative interaction strengths; this is the
principal case of interest in physics.)  The values of the heights
are well-defined up to a global additive constant, so we arbitrarily
pick one of the heights and set its value to be zero.  The heights
$x_i$ and $x_j$ act like they're bound together by a spring with
spring-constant $F_{i,j}$, so that the force pulling $x_i$ and $x_j$
together is $F_{i,j} |x_i-x_j|$, and the energy in the spring is
$\frac12 F_{i,j}(x_i-x_j)^2$.  The total energy of the system is then
  $$E = \sum_{i<j} \frac12 F_{i,j} (x_i-x_j)^2.$$
The probability distribution is (relative to Lebesgue measure) proportional
to $e^{-E}$.

The interaction graph on the sites has an edge between two sites $i$
and $j$ if $F_{i,j}\neq 0$.  We will assume that the interaction graph
contains a spanning tree, since otherwise the system would break apart
into disjoint non-interacting subsystems, which can be dealt with
separately.

The simplest example occurs when $n=2$, where $x_2-x_1$ is distributed
as a normal random variable with mean $0$ and variance $1/F_{1,2}$.

\fref{field} shows a random autonormal / free field configuration where the nonzero springs form a regular 2D grid on the torus.

\begin{figure}[phtb]
\centerline{\epsfig{figure=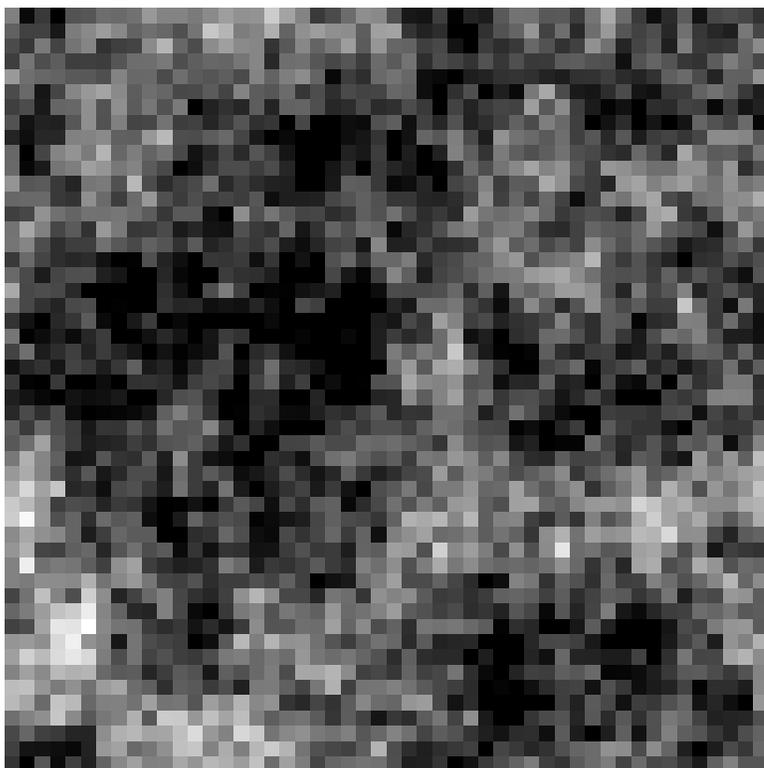,width=4in}}
\caption{
Random free field configuration, with shades of gray representing the
height variables.  The interaction graph is the regular $50\times50$
toroidal grid, and the upper-left-most height is tied to $0$.  This
free field is ``massless''; a massive free field would have an extra
vertex connected to every site on the grid.}
\label{fig:field}
\end{figure}

The reason it's called a {\it free\/}
field (as opposed to another kind of field) is that the springs are
ideal, i.e.\ that the force restoring a value to its mean is {\it linear\/}
in the displacement, without higher order terms.

\subsubsection{Gibbs sampling}

Consider the Gibbs-sampling algorithm (single site heat bath).
When the heights at all sites other than site $i$ are fixed,
the total energy is the following quadratic polynomial in $x_i$:
$$ \sum_j \frac12 F_{i,j} (x_i-x_j)^2 + \sum_{i\neq j<k\neq i} \frac12 F_{j,k}(x_j-x_k)^2.$$
The energy is minized when
 $$\sum_j  F_{i,j} (x_i-x_j) = 0,\text{\ \ \ \ i.e.\ when\ \ } x_i = \sum_j x_j F_{i,j} / \sum_j F_{i,j},$$
and the coefficient of $x_i^2$ is $\frac12
\sum_j F_{i,j}$.
Thus the conditional distribution of the height $x_i$ given the
remaining heights is governed by a normal distribution with mean
$$\sum_j x_j F_{i,j} / \sum_j F_{i,j},$$
and variance $$1/\sum_j F_{i,j}.$$
Equivalently, the height $x_i$ acts as if a spring with spring
constant $\sum_j F_{i,j}$ is pulling it to a weighted average of the
neighboring heights.
Recall that the Gibbs sampler visits the sites, either in sequence or
at random, and randomizes the height $x_i$ at a visited site $i$ by
drawing it from the conditional distribution given the remaining heights.
The term ``autonormal'' comes from that fact that each variable is
normally distributed with nonrandom variance and mean determined by a
weighted average of its neighbors.

\subsubsection{Linear algebra methods}

If the matrix of interaction strengths $F_{i,j}$ can be diagonalized into an orthonormal basis of eigenvectors $\vec{v}_1,\ldots,\vec{v}_n$ with eigenvalues $\lambda_1,\ldots,\lambda_n$, then a random sample can be generated by
$$ \sum_{i:\lambda_i\neq 0} \frac{\Normal(0,1)}{\sqrt{\lambda_i}} \vec{v}_i,$$
where the Gaussian random variables in the sum are independent of one another.  If the interaction graph is a regular lattice, then this approach becomes particularly effective as FFTs can be used (see e.g.\ \citep{dietrich-newsam:gaussian}).  If the interaction strengths are nonuniform or the graph is irregular, then the linear algebra becomes more complicated, and practitioners often prefer the simplicity offered by Markov chain approaches \citep{besag:personal}.

\subsection{Using multishift coupling}

An important property of the layered multishift couplers
is that they are all monotone
couplings.  What this means is that if somehow we can get an upper
bound and lower bound on the values of each variable, $\ell_i\leq x_i\leq u_i$
for each i, then we can repeatedly update the lower configuration
($x_i=\ell_i$ for each $i$) and upper configuration ($x_i=u_i$ for each $i$)
using Gibbsian updates with the
layered multishift coupler for the normal (\sref{normal}).
When the upper and lower configurations get mapped to the same value,
every other possible configuration also gets mapped to this value.
Because we are using a multishift coupler that maps reasonably large
segments of the reals to the same point, the upper and lower
configurations will in fact (with probability $1$) eventually converge
to exactly the same value.  Pseudocode for this approach is given in
\fref{truncated-normal}.

\begin{figure}[phtb]
\begin{code}
\item   $T := 1$\ \ \ \ /* Start at time $-1$ in the past */
\item   Repeat \{
\item\s   /* Set := truncated \sspace */
\item\s   $\ell_1:=0; u_1:=0$\ \ \ \ /* First site is tied to 0 */
\item\s   For $i:=2$ \To $n$
\item\s\s   $\ell_i:=-10^6$; $u_i:= 10^6$\ \ \ \ /* (plausible assumption) */  
\item\s   For $t$ := $T$ DownTo $1$ \{ \ \ \ \ /* Proceed to time zero */
\item\s\s   /* Take care to use previously used random coins */
\item\s\s   If $t$ is a power of 2
\item\s\s\s   SetRandomSeed(seed[$i$,$\log_2(t)$])
\item\s\s   /* ApplyRandomMap(Set) */
\item\s\s   For $i:=2$ To $n$ /* randomize each site (except the one tied to 0) */
\item\s\s\s   /* Apply the multishift coupler for the normal at site $i$ */
\item\s\s\s   /* First pick the parameters defining $f_{L,R,X}()$ at site $i$ */
\item\s\s\s   $\sigma := \sum_j F_{i,j}$
\item\s\s\s   $X := \sigma\Normal(0,1)$
\item\s\s\s   $Y := \exp(-(X/\sigma)^2/2) \Uniform(0,1)$
\item\s\s\s   If ($X<0$) Then $Y := 1-Y$
\item\s\s\s   $L := -\sigma\sqrt{-2\log(1-Y)}$
\item\s\s\s   $R :=  \sigma\sqrt{-2\log(Y)}$
\item\s\s\s   /* Next apply $f_{L,R,X}()$ to the upper and lower bounds */
\item\s\s\s   $\mu_\ell := \sum_j \ell_j F_{i,j} / \sum_j F_{i,j}$
\item\s\s\s   $\mu_u    := \sum_j    u_j F_{i,j} / \sum_j F_{i,j}$
\item\s\s\s   $\displaystyle \ell_i := f_{L,R,X}(\mu_\ell) = \left\lfloor \frac{\mu_\ell+R-X}{R-L}\right\rfloor (R-L) +X$
\item\s\s\s   $\displaystyle  u_i   := f_{L,R,X}(\mu_u)    = \left\lfloor \frac{\mu_u   +R-X}{R-L}\right\rfloor (R-L) +X$
\item\s   \}
\item\s   /* It's now time zero, test for coalescence */
\item\s   If $u_i=\ell_i$ for each $i$ then return $[\ell_1,\ldots,\ell_n]$
\item\s   $T := 2*T$\ \ \ \ /* Otherwise try again starting twice as far in the past */
\item   \}
\end{code}
\caption{Pseudocode for generating random samples from an approximation to the autonormal distribution.  The basic algorithm is monotone coupling from the past using the Gibbs sampler Markov chain, with $\ell_i$ and $u_i$ representing lower and upper bounds on the height at site $i$.  The multishift coupler for the normal distribution is used when doing updates because (1) it correctly give the normal distribution, (2) it is monotone, and (3) with probability one the upper and lower bounds on the state will eventually be exactly equal.  The approximation error (which comes from using $\pm10^6$ in place of $\pm\infty$ in the initial upper and lower bounds on the state) will in practice be minor; in \protect\sref{using-murdoch} we will see how to eliminate this error altogether.}
\label{fig:truncated-normal}
\end{figure}

In the event that some of the interaction strengths $F_{i,j}$ are
negative, we can just use a combination of monotone and anti-monotone
coupling.  In \fref{truncated-normal}, the definitions of $\mu_\ell$
and $\mu_u$ become
$$\mu_\ell := \frac{\sum_{j:F_{i,j}\geq0} \ell_j F_{i,j} + \sum_{j:F_{i,j}<0}    u_j F_{i,j}}{\sum_j F_{i,j}}$$
and
$$\mu_u    := \frac{\sum_{j:F_{i,j}\geq0}    u_j F_{i,j} + \sum_{j:F_{i,j}<0} \ell_j F_{i,j}}{\sum_j F_{i,j}}.$$
\cite{moller:conditional} also mentions mixed monotone/anti-monotone coupling.

How do we get these upper and lower bounds?  It may be tempting to
simply use values such as $10^6$ and $-10^6$, since only a small
portion of the probability distribution is so far out in the tails.
If we did truncate the state space, the running time would increase
logarithmically in the truncation parameter, while the truncation bias
would decrease exponentially.  But there
are reasons not do this: 1) for strongly coupled systems, using such
large values will noticeably and needlessly slow convergence, 2) for
weakly coupled systems $10^6$ may not be large enough, and 3) it's
theoretically displeasing.  In \sref{using-murdoch} we see how to do
without artificial truncations such as this.

\subsection{Using Murdoch's method}
\label{sec:using-murdoch}

At this point the reader should go back and read \sref{murdoch} if
(s)he has not done so already.

\subsubsection{Proposal distribution}
\newcommand{\Et}{E_{\operatorname{tree}}}
There is room for engineering art when picking the proposal
distribution for the independence sampler.
One reasonable choice for the autonormal is
the following.  First pick a spanning tree of the graph rooted at the
special vertex whose height is zero.  We assign the value $x_v$ at vertex
$v$ only after assigning the value at $v$'s parent $u$ in the tree.  The
distribution of $x_v$ is a normal with mean $x_u$ and variance $2/F_{u,v}$.
Call the resulting configuration $B$.  Let $\Et(B)$ be the energy
contained in just those springs that are part of the spanning tree.
The probability density of the proposal distribution (relative to
Lebesgue measure) is then $\exp(-\Et(B)/2)$.

According to the Metropolis-Hastings update rule, when the current state
is $A$ and the proposal is $B$, we always accept the proposal $B$ if
  $$\pi(A) P_A(B) \leq \pi(B) P_B(A),$$
where $\pi(x)$ is the desired probability of state $x$ in a discrete space,
and $P_x(y)$ is the probability of a transition from $x$ to $y$.
Otherwise we accept the proposal with some probability less than one.
In the continuum limit, for our application the above relation amounts to
\begin{align*}
  \exp(-E(A)) \exp(-\Et(B)/2) &\leq \exp(-E(B)) \exp(-\Et(A)/2)\\
              E(A) - \Et(A)/2 &\geq E(B) - \Et(B)/2 \\
        E(A) + (E(A)- \Et(A)) &\geq 2 E(B) - \Et(B) \\
\intertext{which holds whenever}
                         E(A) &\geq E_{\max} \equiv 2 E(B) - \Et(B).
\end{align*}
Any state $A$ with energy $E(A) \geq E_{\max}$ gets mapped to state
$B$.  Furthermore, $E(B) \leq E_{\max}$.  Therefore, after the
Metropolis-Hastings update we are guaranteed that the energy of the
updated state is at most $E_{\max}$.

\subsubsection{Finite box containing updated state}

We again use the spanning tree when converting this bound on the energy
to an upper and lower bound on the value of each coordinate.
Vertices $v$ adjacent to the distinguished vertex have easy bounds on
their values $x_v$:
$$   |x_v| \leq \sqrt{2 E_{\max} / F_{v,0}},$$
since if $|x_v|$ were any larger, the energy in just the spring connecting
$v$ to the distinguished vertex would exceed to total possible energy $E_{\max}$.

To deal with vertices further away from the special vertex, we prove by
induction the following claim:
\begin{claim}
  Given vertices $v_0,\ldots,v_k$, where $k>0$, $x_{v_0}=0$ and
  $x_{v_k} = x$, if we seek to minimize the energy just in the springs
  $(v_0,v_1), \ldots, (v_{k-1},v_k)$, this minimum energy is
     $\frac12 [ 1/(1/F_1 + \cdots + 1/F_k) ] x^2$
  where $F_i$ denotes $F_{v_{i-1},v_i}$.
\end{claim}
\begin{proof}
This claim hold trivially for $k=1$.  Suppose that it holds for $k$, we prove
it for $k+1$.  Let $x=x_{v_{k+1}}$ and $y=x_{v_{k}}$.  By induction
the minimum energy
contained in the given springs is
 $$  \frac12 F y^2 + \frac12 F' (y-x)^2 $$
where we have for convenience let $F$ denote $[ 1/(1/F_1 + ... + 1/F_k) ]$
and $F'$ denote $F_{k+1}$.  This energy is minimized when
\begin{align*}
  F y + F' (y-x) &= 0 \\
               y &= F' x / (F + F')
\end{align*}
at which point the energy takes the value
\begin{align*}
 E &= \frac12 F \left({F' x \over F+F'}\right)^2 +
      \frac12 F' \left({F' x \over F+F'}-x\right)^2           \\
   &= \frac12 F \left({F' x \over F+F'}\right)^2 +
      \frac12 F' \left({-F x \over F+F'}\right)^2                \\
   &= \frac12 \left[{(F {F'}^2 + F^2 F')\over(F+F')^2}\right] x^2 \\
   &= \frac12 \left[{ F F' \over F+F'}\right] x^2               \\
   &= \frac12 \left[{ 1\over 1/F + 1/F'}\right] x^2             \\
   &= \frac12 \left[{ 1\over (1/F_1 + \cdots + 1/F_k) + 1/F_{k+1}} \right] x^2
\end{align*}
as claimed.
\end{proof}

From this we conclude
    $$ |x_{v_k}| \leq \sqrt{2 E_{\max} [1/F_1 + \cdots + 1/F_k]}. $$

\subsection{CFTP using composite random maps}

Next we suitably mix the independence sampler and the Gibbsian updates
to define a composite Markov chain with which we can do CFTP.  We use
a mixing strategy different from the one originally advocated by
\cite{murdoch:exact-bayesian}, since the strategy below is easier to
use for this problem.  The composite Markov chain that we use is given
by the following update rule:

Input: current state $x$
\begin{enumerate}
\item[T1]
   Generate a proposal state for the independence sampler.
\item[T2]
   Ignoring $x$, get upper and lower bounds $u$ and $\ell$ on resulting state
   that would hold regardless of input.
\item[T3]
   Do Gibbsian updates on $u$ and $\ell$ (but not $x$) with the layered multishift coupler until $u=\ell$.
   Let $C$ be the number of Gibbsian updates performed.
\item[B1]
   Generate a proposal state for the independence sampler
   (independent of the previous one).
\item[B2]
   Ignoring $x$, get upper and lower bounds $u$ and $\ell$ on resulting state
   that would hold regardless of input.
\item[MH]
   With the usual Metropolis-Hastings probability,
   either set $x$ to the proposal, or leave it unmodified.
\item[B3]
   Do $C$ Gibbsian updates with the layered multishift coupler starting from the states $u$, $\ell$, and $x$.
\item[R] Output state $x$, and declare coalescence if $u=\ell$.
\end{enumerate}

We make a few observations:
\begin{enumerate}
\item If we pick a random natural number $C$ from any distribution and then do $C$ updates of a state distributed according to $\pi$, the result will be distributed according to $\pi$.  Therefore the above randomizing operation preserves the desired distribution $\pi$.
\item If $u=\ell$, then the output state is independent of the input state.
\item $\Pr[u=\ell] \geq 1/2$.
\end{enumerate}

Now we view these randomizing process as one step of a composite Markov chain with the desired steady-state distribution $\pi$, and do CFTP with the composite randomizing operations.
To do CFTP we compose the maps defined by these Markovian updates
going back in time.  A convenient way to do this is to just keep
trying new random maps $F_{-1}, F_{-2}, \ldots$ until we find a map $F_{-T}$
which by itself is coalescent (as determined by the $u=\ell$ test).
The expected value of $T$ is at most $2$.  Then we take the image of
$F_{-T}$, and determine where the maps
$F_{-T+1},\ldots,F_{-1}$ take it to at time $0$.

We remark that we can also use these composite random maps in the
read-once version of CFTP described by \citep{wilson:rocftp}, and that
this is in fact the approach we took when generating the sample shown
in \fref{field}.

\section*{Source code}

The programs used to make \fref{field} and the simulation results in
\sref{autogamma} are available at this article's web site
\texttt{http://dbwilson.com/shift/}.

\section*{Acknowledgements}

The author thanks Duncan Murdoch for useful comments.  The author also
thanks Jesper M\o ller for calling attention to the autonormal, and
Julan Besag and Tom Spencer for useful background information on it.
The Warwick Randomised Algorithms and Stochastic Simulation (WRASS)
conference and the Fields Institute's workshop on Monte Carlo methods
helped stimulate this article's development.

\nocite{besag-kooperberg:autoregressions}

\bibliography{exact,shift,eus,rocftp}
\bibliographystyle{plainnat}

\end{document}